\documentclass[11pt]{amsart}
\usepackage{graphicx}%
\usepackage{multirow}%
\usepackage{amsmath,amssymb,amsfonts}%
\usepackage{mathrsfs}%
\usepackage{booktabs}%
\usepackage{float}
\usepackage{multirow}
\usepackage{graphicx,epsfig}
\usepackage{color}
\usepackage{enumerate,enumitem}
\usepackage{epsfig,subfigure,epstopdf}
\usepackage{graphicx,euscript}
\usepackage{hyperref} 
\usepackage{multirow,array}
\usepackage[draft,commentmarkup=margin]{changes}
\usepackage{cleveref}

\def\cM{\mathcal{M}}
\def\cN{\mathcal{N}}
\def\vu{\vec{u}}
\def\vnu{\vec{\nu}}
\def\hnu{\hat{\nu}}
\def\vx{\vec{x}}

\def\vX{\vec{X}}
\def\vY{\vec{Y}}
\def\hX{\hat{X}}

\def\vxi{\vec{\xi}}

\def\vtau{\vec{\tau}}
\def\Id{{\rm Id}}
\def\vv{\vec{v}}

\def\vphi{\vec{\varphi}}
\def\half{{1/2}}

\def\vn{\vec{n}}
\def\vw{\vec{w}}
\def\hnu{\hat{{\nu}}}
\def\hn{\hat{n}}

\def\hGamhalf{\hat{\Gamma}_h^{m+\half}}

\def\d{{\rm{d}}}

\def\bfu{\mathbf{u}}
\def\tPsi{\widetilde{\Psi}}
\def\vc{\vec{c}}
\def\ve{\vec{e}}
\def\bfS{\mathbf{S}}
\def\bfW{\mathbf{W}}
\def\bfL{\mathbf{L}}

\definecolor{darkgreen}{rgb}{0,0.5,0}

\definecolor{fixmagenta}{rgb}{0.85,0,0.85}

\theoremstyle{thmstyletwo}%
\newtheorem{theorem}{Theorem}[section]
\newtheorem{proposition}[theorem]{Proposition}%

\newtheorem{remark}{Remark}[section]%
\newtheorem{lemma}{Lemma}[section]
\theoremstyle{definition}
\newtheorem{example}{Example}[section]%
\numberwithin{equation}{section}

\usepackage{cleveref}
\begin{document}


\title[Energy-stable linearized Parametric FEMs]{Surface evolution under conformal map and linearized energy-stable parametric finite element methods}
	\author{Beiping Duan}
	\address{MSU-BIT-SMBU Joint Research Center of  Applied Mathematics, 
	Shenzhen MSU-BIT University, Shenzhen  518172, P.R. China}
\email{duanbeiping@smbu.edu.cn}

\begin{abstract}
In this paper, we aim to develop a linearized, energy-stable, second-order-in-time parametric finite element method without Lagrange multipliers for surface evolution equations on surfaces with or without boundaries, including (volume-preserving) mean curvature flow and surface diffusion. We first reformulate the original equations into a new unified form in which the tangential motion is determined by requiring the flow map to be conformal. Based on this new formulation, we present a family of parametric finite element schemes that employ a linearized Euler or Crank–Nicolson discretization in time. We prove the existence and uniqueness of the resulting linear system, as well as the energy stability of the schemes. Our  discrete methods can be applied to general surfaces although the second-order scheme may degenerate to first-order accuracy in certain cases. Numerical results verify the robustness of our method and demonstrate its ability to preserve mesh quality.
\end{abstract}


\maketitle

\section{Introduction}
Let $\cM$ denote a closed {\it regular surface} in $\mathbb{R}^{3}$ (see \cite[Section 2-2]{do2016differential} for the definition). Assume that $\vX(t): \cM\rightarrow \Gamma(t)$ is a diffeomorphism that parameterizes the closed evolving surface  $\Gamma(t)$ over $\cM$, then the evolution of $\Gamma(t)$ can be described by the flow map $\vX(t)$:
\begin{equation}\label{eqn-original}
	\frac{\partial \vX(t)}{\partial t}\circ\vX^{-1}(t)\cdot \vn=v_n(H),
\end{equation}
where $\vn$ is the unit outward normal vector, and $v_n$ is the normal component of the velocity defined on $\Gamma(t)$.  We  use $v_n(H)$ to emphasize the dependence on the mean curvature $H$. Widely studied geometric flows induced by decreasing the surface area  include:
\begin{itemize}
	\item (volume-preserving) mean curvature flow  $v_n(H)=-H$ ($-H+\langle H\rangle_{\Gamma(t)}$);
	\item surface diffusion \hspace{1.62cm} $v_n(H)=\Delta_{\Gamma(t)}H$;
\end{itemize}
where $\langle \cdot\rangle_{\Gamma(t)}:=\int_{\Gamma(t)}\cdot\,\d_{\Gamma}$. Note that $\int_{\Gamma(t)}Hv_n(H)\, \d_{\Gamma}\le 0$ for all the mentioned cases, so multiplying \eqref{eqn-original} by $H$ and integrating both sides over $\Gamma(t)$ it follows $\frac{d}{dt}|\Gamma(t)|\le 0$, with $|\Gamma(t)|=\int_{\Gamma(t)}1\,\d_{\Gamma}$. 
We refer interested readers to \cite{brakke1978motion,Escher1998,escher1998volume} for analytical results on the flows mentioned above.

\subsection{Background}
Parametric finite element method (FEM), which originates from Dziuk's pioneering work \cite{Dziuk-1990}, has become an important approach in solving geometric evolution equations numerically. For parametric FEMs,, the numerical surface at $t=t_{m+1}$ is determined by the parameterization $\vu_h^{m+1}=\sum_{j}\bfu_j^{m+1}\phi_j$, where $\bfu_j^{m+1}\in \mathbb{R}^3$ and $\{\phi_j\}_j$ are the Lagrange basis functions defined on obtained numerical surface $\Gamma_h^m$. Since the nodal positions $\{\bfu_j^{m+1}\}_j$ determined by corresponding numerical scheme move with time stepping, they may distributed improperly such that the mesh distorts or degenerates at some time while the surface it represents remains smooth. This reminds us the continuous system has to be reformulated in a proper form that can induces artificial tangential motions such that the discretization based on it can distributed the nodal positions in a desired manner. For evolving surfaces, work in this direction includes the Barrett-Garcke-N\"urnberg (BGN) method \cite{BGN2008}, the DeTurck flow method \cite{Elliott-Fritz-2017}, the minimal deformation rate (MDR) method \cite{hu2022evolving}, the minimal deformation (MD) method \cite{duan2024new} and the variant DeTurck flow method \cite{duan2024energy,duan2025mesh}. An alternative approach is to apply mesh-smoothing procedures to maintain mesh quality during the simulation, see e.g., \cite{bansch2005finite}.

In addition to mesh quality, another  point of interest is the design of energy-stable schemes. For geometric flows induced by the surface area, linearized energy-decaying schemes can be found in \cite{BGN2008,Gao-Garcke-Li-Tang-2026}, while nonlinear schemes include \cite{duan2025mesh,GAO2025113718,Garcke-Jiang-2025,duan2026second}. For Willmore flow, we refer interested readers to \cite{Zhao2026}. For surface diffusion with moving boundaries, e.g., the solid-state dewetting model proposed in \cite{Jiang2020Sharp}, there is still no linearized energy-decaying second-order-in-time scheme.


On error estimates for geometric flows of surfaces, related work includes \cite{KLLP2017,Kovacs-Li-Lubich-2019,B2021A,li2025dynamic}, in which Huisken's evolution equations are coupled with the original system to control the mean curvature $H$ and unit normal $\vn$. Most recently, Deckelnick and Styles \cite{deckelnick2026error} proved the optimal $H^1$-error bound for the DeTurck method.  For numerical analysis of curve evolution, interested readers may refer to, e.g., \cite{deckelnick2025finite,deckelnick2025second}.

\subsection{Related work and our contributions}
Conformal mapping plays an important role in both geometric analysis and computational geometry \cite{CCG_2007}. Numerical conformal mapping is a well-established field and has been widely studied (see \cite{trefethen2025numerical} for a brief review). Recently, by requiring the parametrization $\vu_h^{m+1}$ to be conformal, Rao and Ma \cite{ma2026highorder} constructed an energy-stable Runge--Kutta scheme for mean curvature flow and surface diffusion. The scheme can achieve high-order accuracy in time, but it is fully nonlinear. Denote by $\vu^m(t):=\vX(t)\circ\vX^{-1}(t_m)$ the map from $\Gamma^m:=\Gamma(t_m)$ onto $\Gamma(t)$ for $t\ge t_m$, then their formulation at the continuous level reads:
\begin{equation}\label{eqn-Ma-Rao}
	\left\{
	\begin{aligned}
		|\nabla_{\Gamma^m}\vu^m|^2\partial_t\vu^m\cdot\vn&=|\nabla_{\Gamma^m}\vu^m|^2v_n(H),\\
		\Delta_{\Gamma^m}\vu^m&=\frac{1}{2}|\nabla_{\Gamma^m}\vu^m|^2H\vn,
	\end{aligned}
	\right.
\end{equation}
where $v_n=-H$ for mean curvature flow and $v_n=\Delta_{\Gamma(t)}H$ for surface diffusion.  
Note that $\Delta_{\Gamma(t)}H$ is a conformal invariant, i.e.,
\begin{equation*}
	\frac{1}{2}|\nabla_{\Gamma^m}\vu^m|^2(\Delta_{\Gamma(t)}H)\circ\vu^m=\Delta_{\Gamma^m} (H\circ\vu^m),
\end{equation*}
so in the case of surface diffusion, they adopted the following equivalent form of the first equation in \eqref{eqn-Ma-Rao}:
\begin{equation}\label{eqn-Bao}
	\frac{1}{2}|\nabla_{\Gamma^m}\vu^m|^2\partial_t\vu^m\cdot(\vn\circ\vu^m)=\Delta_{\Gamma^m} (H\circ\vu^m).
\end{equation}
Similar idea was adopted in \cite{bao2026arbitrarily} where isoparametric finite elements was utilized to improve the accuracy in space. However, both their methods can be only applied to closed  surfaces of genus zero since those of higher genus or with more complicated topology are not necessarily conformal equivalent, see e.g.,  \cite{imayoshi1992teichmuller}, so the existence may not be guaranteed. Besides, surfaces with boundaries were also not considered.

In this paper, we keep the normal equation \eqref{eqn-original} unchanged and reform the system as follows:
\begin{equation}\label{eqn-relax}
	\left\{
	\begin{aligned}
		(\partial_t\vX)\circ\vX^{-1}\cdot\vn&=v_n(H),\\
		\frac{\omega}{\sqrt{\det J_0^t}}\Delta_{\Gamma^0}\vX+\frac{1-\omega}{\sqrt{\det J_m^t}}(\Delta_{\Gamma^m}\vu^m)\circ\vX^m&=-(H\vn)\circ\vX,
	\end{aligned}
	\right.
\end{equation}
where $\det J_{k}^t$ with $k=0,1,\cdots,m$ represents the Jacobian determinant between $\Gamma^k$ and $\Gamma(t)$, given by $\int_{\Gamma(t)}\d_{\Gamma}=\int_{\Gamma^k}\sqrt{\det J_k^t}\d_{\Gamma^k}$, and  $\omega\in [0,1]$ is a parameter. Note that the conformality is not utilized in the first equation of \eqref{eqn-relax}. Adopting the Crank-Nicolson scheme in time, we obtain a linearized and energy-stable parametric finite element scheme after carefully choosing the predicted intermediate surface $\hGamhalf$. By naturally extending the conformal map to the boundary, we generalize the approach to solving the solid-state dewetting problem proposed in \cite{Jiang2020Sharp}, and the energy stability is preserved, although in this case the energy includes the surface area and the area occupied on the substrate.

Note that at the continuous level, \eqref{eqn-relax} and \eqref{eqn-Ma-Rao} are equivalent under the assumption that $\vX(t)$ is conformal for any $t\in[0,T]$. However, in the discrete setting, our schemes have the following advantages:
\begin{itemize}
	\item Our numerical schemes based on \eqref{eqn-relax} with $\omega=0$ can handle arbitrary surfaces without genus or boundary restrictions, although the second-order scheme may degenerate to first-order accuracy in certain cases;
	\item The term $\frac{\omega}{\sqrt{\det J_0^t}}\Delta_{\Gamma^0}\vX$ can prevent accumulated errors from polluting the conformality in the tangent space, even when $\omega$ is very small.
\end{itemize}

The rest of this paper is organized as follows. We first present some preliminaries in \Cref{sec-pre}. The weak form for the evolution equation \eqref{eqn-original} of closed surfaces and the corresponding parametric finite element schemes are presented in \Cref{sec-cloased-surface}. Numerical schemes for surfaces with boundaries are given in \Cref{sec-boundary}, where the solid-state dewetting model is considered. Numerical results are presented in \Cref{sec-num} to demonstrate the spatial and temporal convergence rates and the robustness of our method.

\section{Preliminaries}\label{sec-pre}
Suppose that $\cM$ is a regular surface embedded in $\mathbb{R}^3$ with local parametrization $\vphi(\theta^1,\theta^2): \mathbb{R}^2\supset U\rightarrow \vphi(U)\subset \cM$ and denote by  $h:=\vphi^*\mathfrak{e}$ the pull-back metric of the Euclidean metric $\mathfrak{e}$, given by $h_{ij}=\partial_i\vphi\cdot\partial_j\vphi$, where we use $\partial_i$ to denote $\partial_{\theta^i}$ for simplicity. For any smooth scalar or vector-valued function $f$ defined on $\cM$, the surface gradient $\nabla_{\cM}$ is given by 
\begin{equation}\label{gradient}
	(\nabla_{\cM}f)\circ\vphi:=\partial_i(f\circ\vphi)h^{ij}\partial_j\vphi^T,
\end{equation} 
where we used the Einstein summation convention.  The surface Laplacian $\Delta_{\cM}$ is defined by
\begin{equation}\label{Laplace}
	(\Delta_{\cM}f)\circ\vphi:=\frac{1}{\sqrt{\det h}}\partial_j\left(\sqrt{\det h}h^{ij}\partial_i f\right),
\end{equation}
where and here after we will  abuse $h$ as $[h_{ij}]$ for simplicity.

\begin{lemma}\label{lemma-conformal}
	Suppose $\cM$ and $\cN$ are two regular surfaces with or without boundaries, and $\Psi$ is a diffeomorphism mapping from $\cM$ onto $\cN$. If $\Psi$ is conformal, then 
	\begin{equation*}
		\Delta_{\cM}\Psi=-\sqrt{\det J}(H\vn)\circ\Psi,
	\end{equation*}
	where $J$ is the Jacobian matrix such that $\int_{\cN}\,\d_\cN=\int_{\cM}\sqrt{\det J}\,\d_\cM$, $H$ and $\vn$ are the mean curvature and unit normal of $\cN$.
\end{lemma}
\begin{proof}
	Denote by $\tPsi=\Psi\circ\vphi$ and $g:=\tPsi^*\mathfrak{e}$ the pull-back metric given by $g_{ij}=\partial_i\tPsi\cdot\partial_j\tPsi$. Note that $\Psi$ is conformal, so it follows $g=e^{\lambda}h$ with $\lambda:=\lambda(\theta^1,\theta^2)$, which gives $g_{ij}=e^{\lambda}h_{ij}$, $g^{ij}=e^{-\lambda}h^{ij}$ and $\sqrt{\det g}=e^{\lambda}\sqrt{\det h}$. 	By definition
	\begin{equation*}
		\begin{aligned}
			(\Delta_{\cM}\Psi)\circ\vphi&=\frac{1}{\sqrt{\det h}}\partial_j\left(\sqrt{\det h}h^{ij}\partial_i\tPsi\right)\\
			&=\frac{e^\lambda}{\sqrt{\det g}}\partial_j\left(\sqrt{\det g}g^{ij}\partial_i\tPsi\right)\\
			&=e^{\lambda}(\Delta_{\cN}\Psi)\circ\vphi=-e^{\lambda}(H\vn)\circ\tPsi.
		\end{aligned}
	\end{equation*}
	The desired equation follows from $e^\lambda=\frac{\sqrt{\det g}}{\sqrt{\det h}}$.
\end{proof}
\begin{lemma}\label{lem-2}
	Let $\cM$, $\cN$ and $\Psi$ be defined as above. Then the following inequality holds
	\begin{equation}\label{eqn-ineqn}
		\frac{1}{2}\int_{\cM}|\nabla_{\cM}\Psi|^2\,\d_{\cM}\ge |\cN|,
	\end{equation}
	where the equal sign holds if and only if $\Psi$ is conformal.
\end{lemma}
\begin{proof}
	For polyhedral surfaces $\cM$ and $\cN$, we refer the reader to \cite[Lemma 2.1]{BGN2008} for a proof. In the general case, it can be easily derived using a local parametrization.  In fact, note that
	\begin{equation*}
		\frac{1}{2}|\nabla_\cM\Psi|^2\circ\vphi =\frac{1}{2}g_{ij}h^{ij}\ge \sqrt{\det (gh^{-1})}=\frac{\sqrt{\det g}}{\sqrt{\det h}}.
	\end{equation*}
	Integrating both sides then yields the desired inequality. The conformality follows from the equal sign holds if and only if $gh^{-1}=e^{\lambda}I$ with $I$ the identity matrix.
\end{proof}

\section{Closed surfaces}\label{sec-cloased-surface}
In this section, we discuss the evolution of closed surfaces. We denote $\Gamma^0$ the initial surface and $\vX(t)$ the diffeomorphism from $\Gamma^0$ onto the unknown $\Gamma(t)$.  We consider the case which satisfies
\begin{equation}\label{eqn-assumption}
	\int_{\Gamma(t)}v_n(H) H\,\d_{\Gamma} \le 0.
\end{equation}
One can easily check \eqref{eqn-assumption} holds for the (volume-preserving) mean curvature flow and surface diffusion.
\subsection{Weak formulation}
Let $0=t_0<t_1<\cdots<t_{M-1}<t_M=T$ be a  uniform partition  of $[0,T]$ with $\tau=t_{m}-t_{m-1}$ the time stepsize. Denote by $\vu^m(t):=\vX(t)\circ\vX^{-1}(t_m)$  the map from $\Gamma^m:=\Gamma(t_m)$ onto $\Gamma(t)$. We fix the normal velocity and force the flow map $\vX(t)$ to be conformal for $t\in[0,T]$, then reform \eqref{eqn-original} with unknowns $\vu^{m}$ and $H$ as follows:
\begin{equation}\label{eqn-relax-1}
	\left\{
	\begin{aligned}
		(\partial_t\vX)\circ\vX^{-1}\cdot\vn&=v_n(H)\\
		\frac{\omega}{\sqrt{\det J_0^t}}\Delta_{\Gamma^0}\vX+\frac{1-\omega}{\sqrt{\det J_m^t}}(\Delta_{\Gamma^m}\vu^{m})\circ\vX^m&=-(H\vn)\circ\vX,
	\end{aligned}
	\right.
\end{equation}
where $\omega\in[0,1]$ is independent of the spatial variable. 

\begin{proposition}\label{prop-1}
	If the solution of \eqref{eqn-relax-1} exists and is smooth enough for $m=0,1,\cdots,M-1$, then the flow map $\vX(t)$ is conformal and $\frac{d}{dt}|\Gamma(t)|\le 0$. 
\end{proposition}
\begin{proof}
	Take $\eta=H$ then $\frac{d}{dt}|\Gamma(t)|$ follows from 
	\begin{equation*}
		\begin{aligned}
			\frac{d}{dt}|\Gamma(t)|&=\int_{\Gamma(t)}\nabla_{\Gamma(t)}\Id:\nabla_{\Gamma(t)} (\partial_t\vX\circ\vX^{-1})\,\d_{\Gamma}\\
			&=\int_{\Gamma(t)}\frac{\partial\vX}{\partial t}\circ\vX^{-1}\cdot H\vn \,\d_{\Gamma}=\int_{\Gamma(t)}v_n(H)H\,\d_{\Gamma}\le 0.
		\end{aligned}
	\end{equation*}
	Test the second equation with $\sqrt{\det J_0^t}\vX(t)$ and integrate on $\Gamma^0$ then we get
	\begin{equation*}
		\begin{aligned}
			\omega &\int_{\Gamma^0}|\nabla_{\Gamma^0}\vX|^2\,\d_{\Gamma^0}+(1-\omega)\int_{\Gamma^m}|\nabla_{\Gamma^m}\vu^{m}|^2\,\d_{\Gamma^m}\\
			&=\int_{\Gamma(t)}-\Delta_{\Gamma(t)}\Id\cdot\Id\,\d_{\Gamma}= \int_{\Gamma(t)} \nabla_{\Gamma(t)}\cdot\Id\,\d_{\Gamma}=2|\Gamma(t)|
		\end{aligned}
	\end{equation*}
	for $m=0,1,\cdots,M-1$, where we utilized $\Delta_{\Gamma(t)}\Id=-H\vn$.  On the other hand, appealing to \Cref{lem-2} we know both $\vX(t)$ and $\vu^{m}$ are conformal if $\omega\in(0,1)$. The case of $\omega=0,1$ follows from that the conformality of $\vX(t)$ and $\vu^{m}$ are equivalent. 
\end{proof}

Since \eqref{eqn-relax-1} implies $\vX(t)$ is conformal, it is equivalent to 
\begin{equation}\label{reform-1}
	\left\{
	\begin{aligned}
		\frac{\partial\vX}{\partial t}\circ\vX^{-1}\cdot \vn&= v_n(H),\\
		\frac{1}{\sqrt{\det J_0^t}}\Delta_{\Gamma^0}\vX\circ\vX^{-1}&=- H\vn.
	\end{aligned}
	\right.
\end{equation}
We would also like to discuss the connection and difference between \eqref{reform-1} and the method proposed by Duan and Li in \cite{duan2024new}, in which the unknowns are  the flow map $\vX(t)$,  the curvature $H$, and the scalar function $\kappa$. The reformulated system in \cite{duan2024new} is  
\begin{equation}\label{duan-Li-2024}
	\left\{
	\begin{aligned}
		\frac{\partial\vX}{\partial t}\circ\vX^{-1}\cdot \vn&= v_n(H),\\
		\Delta_{\Gamma^0}\vX\circ\vX^{-1}&=\kappa \vn,\\
		\Delta_{\Gamma(t)}\Id\cdot \vn& =-H.
	\end{aligned}
	\right.
\end{equation}
The second equation in \eqref{duan-Li-2024} enforces the flow map $\vX(t)$ to be harmonic during the evolution, where $\kappa$ acts as a Langrange multiplier.  Note that a conformal map is always harmonic while the converse does not hold unless the surface is of genus zero. So if the initial surface is of genus zero, the evolution equation \eqref{duan-Li-2024} coincides with \eqref{reform-1} and one can get 
\begin{equation*}
	\kappa=(\Delta_{\Gamma^0}\vX)\circ\vX^{-1}\cdot\vn=-\sqrt{\det J_0^t}H,
\end{equation*}
where the last step follows from \Cref{lemma-conformal}. For the initial surface is not of genus zero, \eqref{duan-Li-2024} and \eqref{reform-1} are different.

Let $\cM$ and $\cN$ be defined as before with $\Psi$ the diffeomorphism mapping form $\cM$ onto $\cN$.  For any $\vY\in\mathbf{H}^1(\cM)^3$ and $\vec{\chi}\in \mathbf{H}^1(\cN)^3$, we introduce the bilinear form 
\begin{equation}\label{we_Conformal}
	\begin{aligned}
		&C_\omega(\vY,\vec{\chi};\cM,\cN) \\
		=&\omega\int_{\cM}\nabla_{\cM}\vY: \nabla_{\cM} (\vec{\chi}\circ\Psi) \,\d_{\cM}+(1-\omega)\int_{\cN}\nabla_{\cN}(\vY\circ\Psi^{-1}):\nabla_{\cN} \vec{\chi}\,\d_{\cN}.
	\end{aligned}
\end{equation}
Our weak form reads: Find $\vX(t)\in \mathbf{H}^1(\Gamma^0)^3$ and $H\in \mathbf{H}^1(\Gamma(t))$ (or $H\in \mathbf{L}^2(\Gamma(t))$) such that
\begin{equation}\label{weak-1}
	\left\{
	\begin{aligned}
		\int_{\Gamma(t)}\frac{\partial\vX}{\partial t}\circ\vX^{-1}\cdot \vn\,\eta \,\d_{\Gamma}-\int_{\Gamma(t)} v_n(H)\eta \,\d_{\Gamma}=0&,\\
		C_\omega(\vX,\vxi\circ\vu^m;\Gamma_0,\Gamma^m) -\int_{\Gamma(t)} (H\vn)\cdot\vxi \,\d_{\Gamma}=0&,
	\end{aligned}
	\right.
\end{equation}
holds $\forall \eta\in \mathbf{H}^1(\Gamma(t))$ (or $H\in \mathbf{L}^2(\Gamma(t))$) and $\vxi\in \mathbf{H}^1(\Gamma(t))^3$.

\subsection{Numerical scheme}
Denote by $\Gamma_h^{m}=\cup_k{T_k^m}$ the polyhedral surface produced by the numerical scheme, where we use $\{p_{k,j}^m\}_{j=1}^3$ (ordered counterclockwise) to denote the vertices of the simplex $T_k^m$ and $\bfu^m\in \mathbb{R}^{3,N_p}$ to represent the collocation of nodal positions. Let $\vX_h^m\in \bfS(\Gamma_h^0)^3$ be the parameterization of $\Gamma_h^m$. For any polyhedral surface $\cN_h$, we define
\begin{equation*}
	\bfS(\cN_h)=\left\{w|\; w \mbox{ is continuous  and piecewise linear on each simplex} \right\},
\end{equation*}
and 
\begin{equation*}
	\bfW(\cN_h)=\left\{w|\;w \mbox{ is constant on each simplex}\right\}.
\end{equation*}
For given polyhedral surface $\cN_h$ and $\forall\phi_h\in \bfS(\cN_h)$, the discrete Laplacian $\Delta_{\cN_h}^h$ is defined as 
\begin{equation*}
	-\left(\Delta_{\cN_h}^h\phi_h,{\chi}_h\right)_{\cN_h}^h=\left(\nabla_{\cN_h}\phi_h,\nabla_{\cN_h}{\chi}_h\right),\quad \forall {\chi}_h\in\bfS(\cN_h),
\end{equation*}
where $(\cdot,\cdot)_{\cN_h}^h$ represents the mass-lumped inner product. Besides, let $\pi_h$ denote the $L^2$-projection from $\bfL^2(\cN_h)^3$ onto $\bfS(\cN_h)^3$, defined by 
\begin{equation*}
	\left(\pi_h\vw, \vxi_h\right)_{\cN_h}^h=\left(\vw, \vxi_h\right)_{\cN_h}, \quad \forall\vxi_h\in \bfS(\cN_h)^3.
\end{equation*}


{\bf Euler scheme}  Find $\vu_h^{m+1}\in \bfS(\Gamma_h^m)^3$ and $H_h^{m+1}\in \bfS(\Gamma_h^m)$ such that
\begin{equation}\label{eqn-Euler}
	\left\{
	\begin{aligned}
		\frac{1}{\tau}\left(\vu_h^{m+1}-\Id,\eta_h\vn_h^{m}\right)^h_{\Gamma_h^m}-\left(v_n(H_h^{m+1}),\eta_h\right)_{\Gamma_h^m}&=0,\\
		C_\omega(\vX_h^{m+1},\vxi_h;\Gamma_h^0,\Gamma_h^{m})-\left( H_h^{m+1},\vn_h^{m}\cdot\vxi_h\right)_{\Gamma_h^m}^h&=0,
	\end{aligned}
	\right.
\end{equation}
holds $\forall\eta_h\in \bfS(\Gamma_h^m)$ and $\vxi_h\in \bfS(\Gamma_h^m)^3$, where $\vn_h^{m}\in \bfW(\Gamma_h^m)^3$ is the unit outward normal of $\Gamma_h^m$, and hereafter we use $\vn_{h,k}^{m}$ to denote the component on the $k$th simplex ${T}_h^{m}$. In the formulation above, 
\begin{equation*}
	v_n(H_h^{m+1})=-H_h^{m+1}, \quad -H_h^{m+1}+\langle H_h^{m+1}\rangle_{\Gamma_h^m}\quad \mbox{or}\quad \Delta_{\Gamma_h^m}^hH_h^{m+1}
\end{equation*}
corresponds to the mean curvature flow, volume-preserving mean curvature flow or surface diffusion.

{\bf Crank-Nicolson scheme}   Let $t_{m+\half}=t_m+\frac{1}{2}\tau$ and $\vX_h^{m+\half}=\frac{1}{2}(\vX_h^{m}+\vX_h^{m+1})$.  Denote by $\hX_h^{m+\half}\in \bfS(\Gamma_h^0)^3$  an appropriate prediction of the exact surface at $t=t_{m+\half}$ whose image is $\hGamhalf:=\cup \hat{T}_{h,k}^{m+\half}$. For simplicity, denote by $\hat{u}_h^{m+\half}:=\hX_h^{m+\half}\circ(\vX_h^{m})^{-1}$. Let   $\delta_t\vu_h^{m+\half}$ be the finite element function in $\bfS(\hGamhalf)^3$ with nodal vector $\delta_t{\bf u}^{m+\half}:=\frac{1}{\tau}({\bf u}^{m+1}-{\bf u}^{m})$. Our Crank-Nicolson scheme reads: Find $\vu_h^{m+1}\in \bfS(\hGamhalf)^3$ and $H_h^{m+\half}\in \bfS(\hGamhalf)$ such that
\begin{equation}\label{eqn-CN}
	\left\{
	\begin{aligned}
		\left(\delta_t\vu_h^{m+\half},\eta_h\hn_h^{m+\half}\right)^h_{\hGamhalf}-\left(v_n(H_h^{m+\half}),\eta_h\right)_{\hGamhalf}&=0,\\
		C_\omega(\vX_h^{m+\half},\vxi_h\circ\hat{u}_h^{m+\half};\Gamma_h^0,\Gamma_h^m)-\left( H_h^{m+\half},\hn_h^{m+\half}\cdot\vxi_h\right)^h_{\hGamhalf}&=0,
	\end{aligned}
	\right.
\end{equation}
holds $\forall\eta_h\in \bfS(\hGamhalf)$ and $\vxi_h\in \bfS(\hGamhalf)^3$, where $\hn_h^{m+\half}\in \bfW(\hGamhalf)^3$ is the unit outward normal of $\hGamhalf$, and hereafter we use $\hn_{h,k}^{m+\half}$ to denote the component on the $k$th simplex $\hat{T}_h^{m+\half}$. In the formulation above, 
\begin{equation*}
	v_n(H_h^{m+\half})=-H_h^{m+\half},\quad  - H_h^{m+\half}+\langle H_h^{m+\half}\rangle_{\hGamhalf}\quad \mbox{or}\quad \Delta_{\hGamhalf}^hH_h^{m+\half},
\end{equation*}
corresponds to the mean curvature flow, volume-preserving mean curvature flow or surface diffusion. 
\begin{remark}
	For our Crank-Nicolson scheme, obtaining $\hGamhalf$ by direct extrapolation from $\vX_h^{m-1}$ and $\vX_h^m$ may cause numerical oscillations. Instead, we obtain $\hGamhalf$ by solving our Euler scheme on $[t_{m-\half},t_{m+\half}]$, using $\vX_h^{m-\half}$ as the initial data and $\Gamma_h^m$ as the spatial domain for $m\ge 1$. At the first step, we apply our Euler scheme to obtain $\Gamma_h^1$. Alternatively, one may run our Euler scheme for the first several time steps if the initial surface is nonsmooth and the time stepsize is not sufficiently small.

	In the case $\omega=0$, our Euler scheme reduces to the well-known BGN scheme, and our Crank-Nicolson scheme becomes
	\begin{equation}\label{eqn-r-CN}
		\left\{
		\begin{aligned}
			\left(\delta_t\vu_h^{m+\half},\eta_h\hn_h^{m+\half}\right)^h_{\hGamhalf}
			&=\left(v_n(H_h^{m+\half}),\eta_h\right)_{\hGamhalf},\\
			\left(\nabla_{\Gamma_h^m}\vu_h^{m+\half},\nabla_{\Gamma^m_h}(\vxi_h\circ\hat{u}_h^{m+\half})\right)_{\Gamma_h^m}
			&=\left(H_h^{m+\half},\hn_h^{m+\half}\cdot\vxi_h\right)^h_{\hGamhalf},
		\end{aligned}
		\right.
	\end{equation}
	with $\vu_h^{m+\half}=\vX_h^{m+\half}\circ(\vX_h^m)^{-1}$, which can be regarded as a first-order discretization of the following system
	\begin{equation}\label{eqn-reform-3}
		\left\{
		\begin{aligned}
			\frac{\partial \vX(t)}{\partial t}\circ\vX^{-1}(t)\cdot \vn&=v_n(H),\\
			\Delta_{\Gamma(t)}\Id&=-H\vn.
		\end{aligned}
		\right.
	\end{equation}
	Therefore, our Crank-Nicolson scheme remains applicable even when the solution of \eqref{eqn-relax-1} does not exist (recall that if the solution exists, then $\vX(t)$ is conformal); in that case, it reduces to a first-order discretization of \eqref{eqn-reform-3}.
	%
	%
\end{remark}

Before presenting the theoretical results, we introduce the weighted outward normal vectors $\vnu_h^m={\pi_h\vn_h^m}$ and $\hnu_h^{m+\half}={\pi_h\hn_h^{m+\half}}$
for $\Gamma_h^m$  and $\hGamhalf$, respectively. It is easy to check at  the $j$th vertex of $\Gamma_h^m$
\begin{equation}\label{def-vnu0}
	\vnu_{h,j}^m:={\vnu}_{h}^m(p_j)=\frac{\sum_{k,p_j\in T_k^m} |T_k^m|\vn_{h,k}^m}{\sum_{k,p_j\in T_k^m} |T_k^m|}.
\end{equation}



\begin{theorem}\label{them-1}
	Assume that 1) $\Gamma_h^m$  (or $\hGamhalf$) contains no degenerated simplices; 2) ${\rm span}\{\vnu_{h,j}^{m}\}_j=\mathbb{R}^3$ (or ${\rm span}\{\hnu_{h,j}^{m+\half}\}_j=\mathbb{R}^3$). Then our Euler (or Crank-Nicolson) scheme has unique solution at $t=t_{m+1}$, and satisfies the following energy-decaying property 
	\begin{equation}\label{eqn-Energy_CN}
		\frac{\omega}{2}\big\|\nabla_{\Gamma^0_h}\vX_h^{m+1}\big\|^2_{\Gamma_h^0}+(1-\omega)|\Gamma_h^{m+1}|\le \frac{\omega}{2}\big\|\nabla_{\Gamma^0_h}\vX_h^{m}\big\|^2_{\Gamma_h^0}+(1-\omega)|\Gamma_h^{m}|.
	\end{equation}
	In addition, \eqref{eqn-Energy_CN} yields $|\Gamma_h^{m+1}|\le |\Gamma_h^0|$.
\end{theorem}
\begin{proof}
	For the uniqueness, we only consider the Eulerian case as the procedure for the Crank-Nicolson scheme is similar. Since the system is linear, we only need to show the uniqueness. To this end, we consider the corresponding homogeneous system with variables $\vY_h\in \bfS(\Gamma_h^m)^3$ and $K_h\in \bfS(\Gamma_h^m)$: 
	\begin{equation*}
		\left\{
		\begin{aligned}
			\left(\vY_h,\eta_h\vn_h^{m}\right)_{\Gamma_h^m}^h-\tau \left(v_n(K_h),\eta_h\right)_{\Gamma_h^m}=0&,\\
			\omega\left(\nabla_{\Gamma^0_h}(\vY_h\circ\vX_h^{m}),\nabla_{\Gamma^0_h}(\vxi_h\circ\vX_h^{m})\right)_{\Gamma_h^0}+(1-\omega)\left(\nabla_{\Gamma_h^m}\vY_h,\nabla_{\Gamma_h^m}\vxi_h \right)_{\Gamma_h^m}&\\
			-\left( K_h,\vn_h^{m}\cdot\vxi_h\right)_{\Gamma_h^m}^h=0&.
		\end{aligned}
		\right.
	\end{equation*}
	Take $\eta_h=K_h$ and  $\vxi_h=\vY_h$ then we get
	\begin{equation*}
		\omega \|\nabla_{\Gamma^0_h}(\vY_h\circ\vX_h^{m})\|_{\Gamma_h^0}^2+(1-\omega) \|\nabla_{\Gamma^m_h}\vY_h\|_{\Gamma_h^m}^2=\tau\left(v_n(K_h),K_h\right)_{\Gamma_h^m}\le 0,
	\end{equation*}
	which implies 
	\begin{equation*}
		\|\nabla_{\Gamma^0_h}(\vY_h\circ\vX_h^{m})\|_{\Gamma_h^0}\quad\mbox{and (or) }\|\nabla_{\Gamma^m_h}\vY_h\|_{\Gamma_h^m} =0,\quad \mbox{and }\left(v_n(K_h),K_h\right)_{\Gamma_h^m}=0.
	\end{equation*}
	Thus $\vY_h$ must be a constant vector, and $K_h\equiv 0$ in the case of mean curvature flow. In the case of surface diffusion or volume-preserving mean curvature flow, $K_h$ must be a constant. Ether case gives $ v_n(K_h)\equiv 0$ and 
	\begin{equation}\label{eqn-Kh}
		\left( K_h,\vn_h^{m}\cdot\vxi_h\right)_{\Gamma_h^m}^h= K_h\left(\vn_h^{m},\vxi_h\right)_{\Gamma_h^m}^h \equiv 0\quad \forall\vxi_h\in \bfS(\Gamma_h^m)^3,
	\end{equation}
	and $ v_n(K_h)\equiv 0$ further implies 
	\begin{equation}\label{eqn-Yh}
		\left( \vY_h,\eta_h\vn_h^{m}\right)_{\Gamma_h^m}^h\equiv 0, \quad \forall \eta_h\in \bfS(\Gamma_h^m).
	\end{equation}
	Choose $\vxi_h=\pi_h\vn_h^m$ and $\eta_h=\phi_j, j=1,2,\cdots, N_p$ with $\phi_j$ the Lagrange basis function at $p_j$. Appealing to assumption 1) and 2) it follows $\vY_h\equiv \vec{0},K_h\equiv 0$ since 
	\begin{equation*}
		\big(\vn_h^{m},\pi_h\vn_h^m\big)_{\Gamma_h^m}^h=\|\vnu_h^m\|_{h,\Gamma_h^m}^2\neq 0,
	\end{equation*}
	and  $\vY_h(p_j)\cdot \vnu_{h,j}$ can not vanish for all $j$, recalling that $\vY_h$ is constant. The uniqueness follows.
	
	Next, we consider the stability of our Crank-Nicolson scheme; the stability of the Euler scheme follows similarly. Choose $\eta_h=H_h^{m+1}$, $\vxi_h=\delta_t\vu_h^{m+1}$ in \eqref{eqn-Euler} and sum the two equations up then it follows
	\begin{equation}\label{eqn-stability-CN}
		\begin{aligned}
			0&\ge \tau\left(v_n(H_h^{m+1}),H_h^{m+1}\right)_{\Gamma_h^m}\\
			&=\frac{\omega }{2} \left(\nabla_{\Gamma^0_h}(\vX_h^{m+1}+\vX_h^m),\nabla_{\Gamma^0_h}(\vX_h^{m+1}-\vX_h^m)\right)_{\Gamma_h^0}\\
			&\quad +\frac{(1-\omega)}{2}\left(\nabla_{\Gamma^m_h}(\vX_h^{m+1}\circ(\vX_h^m)^{-1}+\Id),\nabla_{\Gamma^m_h}(\vX_h^{m+1}\circ(\vX_h^m)^{-1}-\Id)\right)_{\Gamma_h^m}\\
			&= \frac{\omega}{2}\big\|\nabla_{\Gamma^0_h}\vX_h^{m+1}\big\|^2_{\Gamma_h^0}-\frac{\omega}{2}\big\|\nabla_{\Gamma^0_h}\vX_h^{m}\big\|^2_{\Gamma_h^0}\\
			&\qquad +\frac{1-\omega}{2}\|\nabla_{\Gamma^m_h}(\vX_h^{m+1}\circ(\vX_h^m)^{-1})\big\|^2_{\Gamma_h^m}-\frac{1-\omega}{2}\|\nabla_{\Gamma^m_h}\Id\big\|^2_{\Gamma_h^m}.
		\end{aligned}
	\end{equation}
	Then the desired energy inequality  follows from (recalling \eqref{eqn-ineqn})
	\begin{equation*}
		\frac{1}{2}\|\nabla_{\Gamma^m_h}(\vX_h^{m+1}\circ(\vX_h^m)^{-1})\big\|^2_{\Gamma_h^m}\ge |\Gamma_h^{m+1}|, \qquad \frac{1}{2}\|\nabla_{\Gamma^m_h}\Id\big\|^2_{\Gamma_h^m}=|\Gamma_h^m|. 
	\end{equation*}
	The stability of surface area  $|\Gamma_h^{m+1}|\le |\Gamma_h^{0}|$ holds since 
	\begin{equation*}
		\frac{1}{2}\big\|\nabla_{\Gamma^0_h}\vX_h^{m+1}\big\|^2_{\Gamma_h^0}\ge |\Gamma_h^{m+1}|,\qquad \frac{1}{2}\big\|\nabla_{\Gamma^0_h}\Id\big\|^2_{\Gamma_h^0}=|\Gamma_h^0|.
	\end{equation*}

\end{proof}

\section{Surfaces with boundaries}\label{sec-boundary}
For the (volume-preserving) mean curvature flow with fixed boundaries, the weak form has no difference with \eqref{weak-1}. In this section we discuss the solid-state dewetting problem on a plane substrate, which is described by surface diffusion  together with proper boundary conditions \cite{Jiang2020Sharp}. Since boundary is involved, we use $\Gamma(t)$ to represent the interior domain with $\partial\Gamma:=\partial\Gamma(t)$ its boundary, and $\bar{\Gamma}(t)=\Gamma(t)\cup \partial\Gamma$. We still use $\vn$ to denote the unit outward normal of the surface with $\vn|_{\partial\Gamma}$ its restriction on $\partial\Gamma$. Denote by $\vn_{\partial\Gamma}$ (lying the $xoy$-plane) the unit outward normal of $\partial\Gamma$. 

In addition, for simplicity we denote $\vc_{\partial\Gamma}=\vtau_{\partial\Gamma}\times\vn|_{\partial\Gamma}$ with $\vtau_{\partial\Gamma}$ the unit tangent vector of $\partial\Gamma$ such that $\vc_{\partial\Gamma}$ points downward. We further denote by $\vc_{\partial\Gamma^0}=\vc_{\partial\Gamma}|_{t=0}$. This model can be regarded as an $H^{-1}$ gradient flow associated with the energy $|\Gamma(t)|-\cos\alpha|\Omega(t)|$, where $\alpha=\arccos(\vc_{\partial\Gamma}\cdot\vn_{\partial\Gamma})$ is the contact angle, and $\Omega(t)$ is the enclosed area of $\partial\Gamma$.

Similar as before, by forcing $\vX(t)$ to be conformal we reform the evolution equation as follows
\begin{equation}\label{reform-2}
	\left\{
	\begin{aligned}
		\frac{\partial\vX}{\partial t}\circ\vX^{-1}\cdot \vn&= \Delta_{\Gamma(t)}H,\\
		\frac{\omega}{\sqrt{\det J_0^t}}\Delta_{\Gamma^0}\vX+\frac{1-\omega}{\sqrt{\det J_m^t}}(\Delta_{\Gamma^m}\vu^{m})\circ\vX^m&=-(H\vn)\circ\vX,
	\end{aligned}
	\right.
\end{equation}
where all the notations share the same meaning as before. 

The  boundary (contact line) $\partial{\Gamma}$ satisfies:
\begin{enumerate}[label=(\roman*)]
	\item[i)] The contact line condition: $\vX(t)|_{\partial\Gamma}$ is on  the $xoy$ plane and the $z$-component of $\vX(t)$ is $\ge 0$;
	\item[ii)] The relaxed contact angle condition
	\begin{equation*}
		\partial_t\vX|_{\partial\Gamma}\cdot\vn_{\partial\Gamma}=-\beta(\vc_{\partial\Gamma}\cdot\vn_{\partial\Gamma}-\cos\alpha),\quad \beta>0;
	\end{equation*} 
	\item[iii)] The zero-mass flux condition
	\begin{equation*}
		(\vc_{\partial\Gamma}\cdot\nabla_{\Gamma(t)}H)|_{\partial\Gamma}=0.
	\end{equation*}
\end{enumerate}
The tangent motion on the boundary will be determined by the natural continuation of the conformal map from the interior domain to the boundary.
\subsection{Weak formulation}
We next present the weak formulation based on \eqref{reform-2}. Denote by
\begin{equation*}
	\mathbf{H}^1_{z}(\Gamma(t))^3=\mathbf{H}^1(\Gamma(t))\times \mathbf{H}^1(\Gamma(t))\times \mathbf{H}_0^1(\Gamma(t)).
\end{equation*}
For any $\vxi\in \mathbf{H}^1_{z}(\Gamma(t))^3$ we have
\begin{equation}\label{eqn-weak-c}
	\begin{aligned}
		&\int_{\Gamma(t)} \frac{1}{\sqrt{\det J_0^t}}(\Delta_{\Gamma^0}\vX)\circ\vX^{-1}\cdot \vxi\,\d_{\Gamma}\\
		&=\int_{\Gamma^0} \Delta_{\Gamma^0}\vX\cdot (\vxi\circ\vX)\,\d_{\Gamma^0}\\
		&=-\int_{\Gamma^0} \nabla_{\Gamma^0}\vX: \nabla_{\Gamma^0}(\vxi\circ\vX)\,\d_{\Gamma^0}+\int_{\partial\Gamma^0}\nabla_{\Gamma^0}\vX\vc_{\partial\Gamma^0}\cdot (\vxi\circ\vX)\,\d_{\partial\Gamma^0}.
	\end{aligned}
\end{equation}
and similarly
\begin{equation}\label{eqn-weak-id}
	\begin{aligned}
		&\int_{\Gamma(t)} \frac{1}{\sqrt{\det J_m^t}}(\Delta_{\Gamma^m}\vu^{m})\circ(\vu^{m})^{-1}\cdot \vxi\,\d_{\Gamma}\\
		&=-\int_{\Gamma^m} \nabla_{\Gamma^m}\vu^{m}: \nabla_{\Gamma^m}(\vxi\circ\vu^{m})\,\d_{\Gamma^m}+\int_{\partial\Gamma^m}\nabla_{\Gamma^m}\vu^{m}\vc_{\partial\Gamma^m}\cdot (\vxi\circ\vu^{m})\,\d_{\partial\Gamma^m}.
	\end{aligned}
\end{equation}
\begin{lemma}\label{boundary-term}
	Suppose $\vX(t)$ is conformal from $\Gamma^0$ onto $\Gamma(t)$ for $t\in[0,T]$ and naturally extended to the boundary, then it follows
	\begin{equation*}
		\int_{\partial\Gamma^0}(\nabla_{\Gamma^0}\vX\vc_{\partial\Gamma^0})\cdot \vxi\circ\vX\,\d_{\partial\Gamma^0}=\int_{\partial\Gamma}\vc_{\partial\Gamma}\cdot \vxi\,\d_{\partial\Gamma}.
	\end{equation*}
	The same result holds for the boundary term  in \eqref{eqn-weak-id}.
\end{lemma}
\begin{proof}
	The operator $\nabla_{\Gamma^0}\vX$ push any $\vw\in T_{p}\Gamma^0$ forward to $T_{\vX(t,p)}\Gamma(t)$. In fact, suppose $\Gamma^0$ is locally  parameterized by $\vphi(\theta^1,\theta^2)$ with $h_{ij}=\partial_i\vphi\cdot\partial_j\vphi$, then
	\begin{equation*}
		(\nabla_{\Gamma^0}\vX)\circ\vphi\partial_k\vphi=h^{ij}\partial_i(\vX\circ\vphi)\partial_j\vphi^T\partial_k\vphi=h^{ij}\partial_i(\vX\circ\vphi)h_{jk}=\partial_k(\vX\circ\vphi).
	\end{equation*}
	So suppose $\vc_{\partial\Gamma^0}=a^k\partial_k\vphi$, then we have
	\begin{equation}\label{eqn1-lemma2}
		(\nabla_{\Gamma^0}\vX\,\vc_{\partial\Gamma^0})\circ\vphi=a^k\partial_k(\vX\circ\vphi).
	\end{equation}
	We suppose $\partial\Gamma^0$ is locally parameterized by $\vphi(s):=\vphi(\theta^1(s),\theta^2(s))$, Then $\partial\Gamma$ can be expressed locally by $\vX\circ\vphi(s)$. Denote by $b_l=\partial_s\theta^l(s)$, recalling that $\vc_{\partial\Gamma^0}\circ\vphi\perp\partial_s\vphi(s)$, then we have  
	\begin{equation}\label{eqn2-lemma2}
		(a^k\partial_k\vphi)\cdot(\partial_l\vphi b^l)=a^kh_{kl}b^l=0.
	\end{equation}
	Since $\vX(t)$ is conformally extended to the boundary, it holds $g_{kl}=e^{\lambda}h_{kl}$ on the boundary which gives $a^kg_{kl}b^l=0$, where $g_{kl}=\partial_k(\vX\circ\vphi)\cdot \partial_l(\vX\circ\vphi)$ . Thus we get 
	\begin{equation}\label{eqn3-lemma2}
		(\nabla_{\Gamma^0}\vX\,\vc_{\partial\Gamma^0})\circ\vphi\perp b^l\partial_l\vX\circ\vphi=\sqrt{b^lb^r g_{lr}}\,\vtau_{\partial\Gamma}\circ\vX\circ\vphi,
	\end{equation}
	where  $\vtau_{\partial\Gamma}$ is the unit tangent of $\partial\Gamma$. Note that  $\vX(t)$ does not change the orientation of $\Gamma(t)$, so
	\begin{equation}\label{eqn4-lemma2}
		\begin{aligned}
			\nabla_{\Gamma^0}\vX\,\vc_{\partial\Gamma^0}&=|\nabla_{\Gamma^0}\vX\,\vc_{\partial\Gamma^0}|\vc_{\partial\Gamma}=\sqrt{a^ka^rg_{kr}}\vc_{\partial\Gamma}\\
			&=e^{\lambda/2}\sqrt{a^ka^rh_{kr}}\vc_{\partial\Gamma}=e^{\lambda/2}\vc_{\partial\Gamma},
		\end{aligned}
	\end{equation}
	where we utilized $|\vc_{\partial\Gamma^0}|=|\vc_{\partial\Gamma}|=1$ and $|\vc_{\partial\Gamma^0}|=\sqrt{a^ka^rh_{kr}}$. Notice that
	\begin{equation*}
		\frac{|\partial_s(\vX\circ\vphi(s))|^2}{|\partial_s\vphi(s)|^2}=\frac{b^lb^rg_{lr}}{b^lb^rh_{lr}}=e^{\lambda},
	\end{equation*}
	then it follows from \eqref{eqn4-lemma2}
	\begin{equation}\label{eqn5-lemma2}
		\nabla_{\Gamma^0}\vX\,\vc_{\partial\Gamma^0}=\frac{|\partial_s(\vX\circ\vphi(s))|}{|\partial_s\vphi(s)|}\vc_{\partial\Gamma}\circ\vX,
	\end{equation}
	from which we can get the desired equality.
\end{proof}

Test the first equation of \eqref{reform-2} by $\eta\in \mathbf{H}^1(\Gamma(t))$ and the second one by $\vxi\in \mathbf{H}^1_{z}(\Gamma(t))^3$ with using the boundary condition iii) and Lemma \ref{lemma-conformal} we get the weak form: Find $\vX(t)\in \mathbf{H}^1_{z}(\Gamma(0))^3$ and $H\in\mathbf{H}^1(\Gamma(t))$ such that 
\begin{equation}\label{weak-2}
	\left\{
	\begin{aligned}
		\int_{\Gamma(t)}\frac{\partial\vX}{\partial t}\circ\vX^{-1}\cdot \vn\,\eta \,\d_{\Gamma}+\int_{\Gamma(t)} \nabla_{\Gamma(t)} H \cdot \nabla_{\Gamma(t)}\eta \,\d_{\Gamma}=&0,\\
		C_\omega(\vX,\vxi\circ\vu^m;\Gamma^0,\Gamma^m)-\int_{\Gamma(t)} (H\vn)\cdot\vxi \,\d_{\Gamma}=&\\
		\cos\alpha\int_{\partial\Gamma^0}\vxi\circ\vX\cdot(\partial_s\vX)^\perp \d_{\partial\Gamma^0}-\frac{1}{\beta}\int_{\partial\Gamma} (\vX_t\circ\vX^{-1}\cdot\vn_{\partial\Gamma})(\vxi\cdot\vn_{\partial\Gamma})\,\d_{\partial\Gamma},&
	\end{aligned}
	\right.
\end{equation}
holds $\forall \eta\in \mathbf{H}^1(\Gamma(t))$ and $\vxi\in \mathbf{H}^1_{z}(\Gamma(t))^3$, where the bilinear form $C_\omega(\vX,\vxi\circ\vu^m;\Gamma^0,\Gamma^m)$ is given in \eqref{we_Conformal}. The last two terms in \eqref{weak-2} follows from  applying the same procedures as in \cite{bao2023}. In fact, denote by $\ve_z=[0,0,1]^T$, then for $\vxi\in \mathbf{H}^1_{z}(\Gamma(t))^3$ it holds
\begin{equation}\label{eqn-boundary}
	\begin{aligned}
		&\int_{\partial\Gamma} \vc_{\partial\Gamma}\cdot\vxi\,\d_{\partial\Gamma}\\
		&=\int_{\partial\Gamma}\big((\vc_{\partial\Gamma}\cdot\ve_z)\ve_z+(\vc_{\partial\Gamma}\cdot\vn_{\partial\Gamma})\vn_{\partial\Gamma}\big)\cdot\vxi\,\d_{\partial\Gamma}\\
		&=\int_{\partial\Gamma} (\vc_{\partial\Gamma}\cdot\vn_{\partial\Gamma})(\vxi\cdot\vn_{\partial\Gamma})\,\d_{\partial\Gamma}\\
		&=\int_{\partial\Gamma} (\cos\alpha-\beta^{-1}\vX_t\circ\vX^{-1}\cdot\vn_{\partial\Gamma})(\vxi\cdot\vn_{\partial\Gamma})\,\d_{\partial\Gamma}\\
		&=\cos\alpha\int_{\partial\Gamma^0}\vxi\circ\vX\cdot(\partial_s\vX)^\perp \d_{\partial\Gamma^0}-\frac{1}{\beta}\int_{\partial\Gamma} (\vX_t\circ\vX^{-1}\cdot\vn_{\partial\Gamma})(\vxi\cdot\vn_{\partial\Gamma})\,\d_{\partial\Gamma},
	\end{aligned}
\end{equation}
where we use $s$ to denote the arclength parameter of $\Gamma^0$, and $(\cdot)^\perp$ to represent the anticlockwise rotation.

\begin{proposition}
	Suppose $(\vX(t),H)$ solves \eqref{weak-2} up to $t\le T$ then the following energy-decaying property holds
	\begin{equation}\label{d-Energy}
		\frac{d}{dt}|\Gamma(t)|-\cos\alpha\frac{d}{dt}|\Omega(t)| \le 0.
	\end{equation}
	If the boundary condition ii) is also satisfied, then $\vX(t)$ is conformal. 
\end{proposition}
\begin{proof}
	Take $\eta=H$ and $\vxi=\partial_t\vX(t)\circ\vX^{-1}$, add the two equations in \eqref{weak-2} up then we get
	\begin{equation}
		\frac{d}{dt}|\Gamma(t)|-\cos\alpha\frac{d}{dt}|\Omega(t)|=-\big\|\nabla_{\Gamma(t)}H\big\|^2_{\Gamma(t)}-\beta^{-1}\big\|\partial_t\vX(t)\circ\vX^{-1}\cdot\vn_{\partial\Gamma}\big\|^2_{\partial\Gamma}\le 0,
	\end{equation} 
	where we utilized 
	\begin{equation*}
		\int_{\partial\Gamma^0}\vX_t(t)\cdot(\partial_s\vX)^\perp \d_{\partial\Gamma^0} =\int_{\partial\Gamma}\vX_t(t)\circ\vX^{-1}\cdot\vn_{\partial\Gamma}\, \d_{\partial\Gamma}=\frac{d}{dt}|\Omega(t)|.
	\end{equation*}
	Take $\vxi=\Id|_{\Gamma(t)}$  with noting that
	\begin{equation*}
		\begin{aligned}
			\int_{\Gamma(t)} H\vn\cdot\vxi\d_{\Gamma}&=-\int_{\Gamma(t)} \Delta_{\Gamma(t)}\Id \cdot \vxi\,\d_{\Gamma}\\
			&=\int_{\Gamma(t)}{\nabla_{\Gamma(t)}\cdot \vxi}\,\d_{\Gamma}-\int_{\partial\Gamma}\vc_{\partial\Gamma}\cdot\vxi\,\d_{\partial\Gamma},
		\end{aligned}
	\end{equation*}
	we then get
	\begin{equation*}
		C_{\omega}(\vX,\vu^m;\Gamma^0,\Gamma^m)=\int_{\Gamma(t)}{\nabla_{\Gamma(t)}\cdot \Id}\,\d_{\Gamma}=2|\Gamma|,
	\end{equation*}
	where the boundary terms are canceled by appealing to \eqref{eqn-boundary} and the boundary condition ii). Recalling \eqref{lemma-conformal} we get the conformality of $\vX(t)$.
\end{proof}

\subsection{Numerical scheme}
Next, we present the numerical schemes based on \eqref{weak-2}. For any polygonal surface $\cN_h$ with nonnegative $z$-component and plane boundary curve $\partial\cN_h$ located in the $xoy$ plane, we use $\vn_{\partial\cN_h}$ to represent the unit outward normal of $\partial\cN_h$ and $\vtau_{\partial\cN_h}$ to denote the unit tangent such that $\vc_{\partial\cN_h}=\vtau_{\partial\cN_h}\times\vn_{\partial\cN_h}$ points downward. Let $\bfS_0(\cN_h)=\bfS(\cN_h)\cap\mathbf{H}_0^1(\cN_h)$ and 
\begin{equation*}
	\bfS_{z}(\cN_h)^3=\bfS(\cN_h)\times\bfS(\cN_h)\times \bfS_0(\cN_h).
\end{equation*}

{\bf Euler scheme} Find $\vu_h^{m+1}\in \bfS_{z}(\Gamma_h^{m})^3$ and $H_h^{m+1}\in \bfS(\Gamma_h^m)$ such that
\begin{equation}\label{Euler-SD}
	\left\{
	\begin{aligned}
		\left(\delta_t\vu_h^{m+1},\eta_h\vn_h^{m}\right)^h_{\Gamma_h^m}+\left(\nabla_{\Gamma_h^m}H_h^{m+1},\nabla_{\Gamma_h^m}\eta_h\right)_{\Gamma_h^m}=&0,\\
		C_\omega(\vX_h^{m+1},\vxi_h;\Gamma_h^0,\Gamma_h^m)-\left( H_h^{m+1},\vn_h^{m}\cdot\vxi_h\right)^h_{_{\Gamma_h^m}}=&\\
		{\cos\alpha}\left(\Big(\partial_s\vX_h^{m+\half}\Big)^\perp,\vxi_h\right)_{\partial\Gamma_h^0}-\frac{1}{\beta}\left(\delta_t\vu_h^{m+1}\cdot\vn_{\partial\Gamma_h^m},\vxi_h\cdot\vn_{\partial\Gamma_h^m}\right)_{\partial\Gamma_h^m},
	\end{aligned}
	\right.
\end{equation}
holds $\forall\vxi_h\in \bfS_z(\Gamma_h^m)^3$ and $\forall\eta_h\in \bfS(\Gamma_h^m)$.

{\bf  Crank-Nicolson scheme} Find $\vu_h^{m+1}\in \bfS_{z}(\hGamhalf)^3$ and $H_h^{m+\half}\in \bfS(\hGamhalf)$ such that
\begin{equation}\label{CN-SD}
	\left\{
	\begin{aligned}
		\left(\delta_t\vu_h^{m+\half},\eta_h\hn_h^{m+\half}\right)^h_{\hGamhalf}+\left(\nabla_{\hGamhalf}H_h^{m+\half},\nabla_{\hGamhalf}\eta_h\right)_{\hGamhalf}=&0,\\
		C_\omega(\vX_h^{m+\half},\vxi_h\circ\hat{u}_h^{m+\half};\Gamma_h^0,\Gamma_h^m)-\left( H_h^{m+\half},\hn_h^{m+\half}\cdot\vxi_h\right)^h_{_{\hGamhalf}}=&\\
		\cos\alpha\left(\Big(\partial_s\vX_h^{m+\half}\Big)^\perp,\vxi_h\right)_{\partial\Gamma_h^0}-\frac{1}{\beta}\left(\delta_t\vu_h^{m+1}\cdot\vn_{\partial\hGamhalf},\vxi_h\cdot\vn_{\partial\hGamhalf}\right)_{\partial\hGamhalf},
	\end{aligned}
	\right.
\end{equation}
holds $\forall\vxi_h\in \bfS_z(\hGamhalf)^3$ and $\eta_h\in \bfS(\hGamhalf)$. In the both schemes above we adopt the same notations as in the previous section.
\begin{remark}
	In our Euler scheme, we use $\partial_s\vX_h^{m+\half}$ in the boundary term rather than $\partial_s\vX_h^{m+1}$ to get the energy-decaying property in \Cref{them-2}.
\end{remark}
\begin{theorem}\label{them-2}
	Assume that 1) $\Gamma_h^m$  (or $\hGamhalf$) contains no degenerated simplices; 2) ${\rm span}\{\vnu_{h,j}^{m}\}_j=\mathbb{R}^3$ (or ${\rm span}\{\hnu_{h,j}^{m+\half}\}_j=\mathbb{R}^3$). Then there exists some  $\epsilon>0$ such that for $\alpha\in (\frac{\pi}{2}-\epsilon,\frac{\pi}{2}+\epsilon)$ the Euler scheme (or the Crank-Nicolson scheme) yields a unique solution and  satisfies
	\begin{equation*}
		\begin{aligned}
			&\frac{\omega}{2}\big\|\nabla_{\Gamma^0_h}\vX_h^{m+1}\big\|^2+(1-\omega) |\Gamma_h^{m+1}| -\cos\alpha|\Omega_h^{m+1}|\\
			\le&\frac{\omega}{2} \big\|\nabla_{\Gamma^0_h}\vX_h^{m}\big\|^2+(1-\omega) |\Gamma_h^{m}| -\cos\alpha|\Omega_h^{m}|.
		\end{aligned}
	\end{equation*}
\end{theorem}
\begin{proof}
	Since the both schemes are linear, the existence and uniqueness follow from considering the corresponding homogeneous system. For our Euler scheme, it reads find $\vY_h\in \bfS_z(\Gamma_h^m)^3$ and $K_h\in\bfS(\Gamma_h^m)$ such that 
	\begin{equation*}
		\left\{
		\begin{aligned}
			\left(\vY_h,\eta_h\vn_h^{m}\right)_{\Gamma_h^m}^h+\tau \left(\nabla_{\Gamma_h^m}K_h,\nabla_{\Gamma_h^m}\eta_h\right)_{\Gamma_h^m}=0&,\\
			\omega\left(\nabla_{\Gamma^0_h}(\vY_h\circ\vX_h^{m}),\nabla_{\Gamma^0_h}(\vxi_h\circ\vX_h^{m})\right)_{\Gamma_h^0}+(1-\omega)\left(\nabla_{\Gamma_h^m}\vY_h,\nabla_{\Gamma_h^m}\vxi_h \right)_{\Gamma_h^m}&\\
			-\left( K_h,\vn_h^{m}\cdot\vxi_h\right)_{\Gamma_h^m}^h +\frac{1}{\beta}\left(\vY_h\cdot\vn_{\partial\Gamma_h^m},\vxi_h\cdot\vn_{\partial\Gamma_h^m}\right)_{\partial\Gamma_h^m}&\\
			-{\cos\alpha}\left(\partial_s(\vY_h\circ\vX_h^m)^\perp,\vxi_h\circ\vX_h^m\right)_{\partial\Gamma_h^0}=0&.
		\end{aligned}
		\right.
	\end{equation*}
	Take $\eta_h=K_h$ and $\vxi_h=\vY_h$, and sum the two equations up then we have for $\alpha=\frac{\pi}{2}$
	\begin{equation*}
		\begin{aligned}
			&\omega \big\|\nabla_{\Gamma^0_h}(\vY_h\circ\vX_h^{m})\big\|_{\Gamma_h^0}^2+(1-\omega) \big\|\nabla_{\Gamma^m_h}\vY_h\big\|_{\Gamma_h^m}^2+\beta^{-1}\big\|\vY_h\cdot\vn_{\partial\Gamma_h^m}\big\|_{\partial\Gamma_h^m}^2\\
			=&-\tau\left\| \nabla_{\Gamma_h^m}K_h\right\|^2_{\Gamma_h^m}\le 0,
		\end{aligned}
	\end{equation*}
	which implies both sides should vanish. So we get $\vY_h$ and $K_h$ are constant. Following similar procedure as we did for  closed surfaces, we get the uniqueness for $\alpha=\frac{\pi}{2}$. Appealing to the perturbation theory for linear systems we know the uniqueness still holds for $\alpha$ in some neighbourhood of $\frac{\pi}{2}$.
	
	Next, we consider the stability of our Crank-Nicolson scheme; the Euler scheme is treated similarly. We take $\eta_h=H_h^{m+\half}$ and $\vxi_h=\delta_t\vu_h^{m+\half}$. We focus only on the boundary terms, as the other terms are handled in the same way as for closed surfaces. For $\zeta\in [t_m,t_{m+1}]$, denote by $\Gamma_h^{\zeta}$ the image of 
	\begin{equation*}
		\frac{t_{m+1}-\zeta}{\tau}\vX_h^m+\frac{\zeta-t_m}{\tau}\vX_h^{m+1}:=\vX_h^\zeta
	\end{equation*}
	with $\vn_{\partial\Gamma_h^\zeta}$ the unit outward normal of the boundary ${\partial\Gamma_h^\zeta}$. We further use $\Omega_h^\zeta$ to denote the enclosed area by the boundary and let $\vv_h^\zeta=\partial_\zeta\vX_h^\zeta=\frac{1}{\tau}(\vX_h^{m+1}-\vX_h^m)$. 
	Note that
	\begin{equation*}
		\begin{aligned}
			|\Omega_h^{m+1}|-|\Omega_h^m|&=\int_{t_{m}}^{t_{m+1}}\frac{d}{d\zeta}|\Omega_h^\zeta|d\zeta\\
			&=\int_{t_{m}}^{t_{m+1}}\int_{\partial\Gamma_h^\zeta}\vv_h^\zeta\circ(\vX_h^{\zeta})^{-1} \cdot \vn_{\partial\Gamma_h^\zeta}\,\d_{\partial\Gamma_h^\zeta} d\zeta\\
			&=\int_{\partial\Gamma_h^0}\int_{t_{m}}^{t_{m+1}}\vv_h^\zeta \cdot {(\partial_s\vX_h^\zeta)^\perp}\,d\zeta\d_{\partial\Gamma_h^\zeta} \\
			&=\int_{\partial\Gamma_h^0}\vv_h^\zeta \cdot \frac{{\tau(\partial_s\vX_h^m+\partial_s\vX_h^{m+1})^\perp}}{2}\,\d_{\partial\Gamma_h^0}\\
			&=\tau\int_{\partial\Gamma_h^0}(\partial_s\vX_h^{m+\half})^\perp\cdot \delta_t\vu_h^{m+\half}\,\d_{\partial\Gamma_h^0},
		\end{aligned}
	\end{equation*}
	where $s$ is the arclength parameter of $\partial\Gamma_h^0$ and in the last step we utilized that the trapezoidal rule exactly expresses the integral when the integrand is linear. The stability follows since the term involving $\beta$ is positive in ether Euler or Crank-Nicolson scheme.
\end{proof}
\section{Numerical results}\label{sec-num}
In this section, we present some numerical experiments for different geometric flows with or without boundaries, to demonstrate the energy-decaying properties of our method, as well as the advantage of mesh distribution. Throughout this section, $N_p$ and $N_T$ denote the numbers of vertices and plane triangles in the mesh, respectively. Moreover, let  $U(t_m)$ denote the numerical solution at time $t=t_m$, with the $j$th component $U_j(t_m)\in \mathbb{R}^{3}$ representing the position of the $j$th vertex. The mesh quality of a two-dimensional polyhedral surface $\cN_h$ is measured by (see e.g., \cite[eq. (7.3)]{Elliott-Fritz-2017})
\begin{equation}\label{msh_q}
	\sigma_{\max}(\cN_h)=\max_{T\in \cN_h}\frac{\mbox{diam } B^*(T)}{\mbox{diam } B_*(T)},
\end{equation}
where $B^*(T)$ and $B_*(T)$ denote the  circumcircle and the maximal inscribed circle of the triangle $T$, respectively. 

For two polyhedral surface $\cM_h$ and $\cN_h$, following \cite[(4.5)]{Bao-Zhao-2021} we define the Hausdorff distance by
\begin{equation}\label{Hausdorff-d}
	\mathcal{H}[\cM_h,\cN_h]=\frac{1}{2}\left[\max_{p\in\cM_h}{\rm dist}(p,\cN_h)+\max_{q\in\cN_h}{\rm dist}(q,\cM_h)\right],
\end{equation}
where ${\rm dist}(p,\cN_h)=\min_{p'\in \cN_h}|p-p'|$ is the distance function, which can be calculated by the algorithm \texttt{point2trimesh} develop by {\it Daniel Frisch} \cite{frisch2016point2trimesh}.

\subsection{Mean curvature flow}
\begin{example}\label{ex1}
	As the first example we test the convergence rates in both space and time. We first consider the standard example of the mean curvature flow, in which the unit sphere $\mathbb{S}^2=\{\vx|\,|\vx|=1\}$ is chosen as the initial surface, then the evolution of $\Gamma(t)$ for $t<0.25$ is given by $\vx(t,\cdot)=\sqrt{1-4t}\Id_{\mathbb{S}^2}$. Denote by $h_0$ the maximal element diameter of $\Gamma_h^0$ and $\tau$ the time stepsize, we define
	\begin{equation*}
		{\rm Error}(t;h_0,\tau)=\max_{j=1,2,\cdots,N_p}\left|U_j(t)-U_j^l(t)\right|
	\end{equation*} 
	with $U_j^l(t)$ the lift of $U_j(t)$ onto $\vx(t,\cdot)$. 
	To test the convergence rate in space, we run our Crank-Nicolson  scheme under different spatial meshes with a uniform and small enough $\tau$ then  compute the convergence rate by 
	\begin{equation}\label{rate_s}
		\mbox{Conver. rate}=\frac{\ln\left({\rm Error}(t;h_0,\tau)/{\rm Error}(t;h_0',\tau)\right)}{\ln(h_0/h_0')},
	\end{equation}
	where we set $t=0.24$. \Cref{CT-MCF-h} presents ${\rm Error}(0.24;h_0,\tau)$ with respect to $h_0$. One can observe that different $\omega$ produces almost the same results with an $\mathcal{O}(h_0^{2})$ convergence rate.
	
	\begin{table}[htp!]
		\centering
		\footnotesize
		\caption{\Cref{ex1}: ${{\rm Error}(0.24;h_0,\tau)}$ under different spatial meshes.}\label{CT-MCF-h}
		\begin{center}
			\begin{tabular}{cccccc}
				\toprule
				$\omega$ & $h_0$   		   & 0.4478    &   0.2273     &  0.1141     &   0.0571		\\
				\toprule
				\multirow{2}{*}{1}        & Error          &1.0370e-01 &  2.0110e-02  & 4.7341e-03  & 1.0751e-03    \\
				& Conver. rate  & --         & 2.42     &  2.10       &  2.14       \\
				\toprule
				\multirow{2}{*}{0.5}    & Error          &1.0362e-01 &  2.0105e-02  & 4.7333e-03  & 1.0748e-03    \\
				& Conver. rate  & --         & 2.42     &  2.10       &  2.14       \\
				\toprule
				\multirow{2}{*}{0.1}     & Error          &1.0320e-01  & 2.0075e-02   & 4.7308e-03  & 1.0745e-03    \\
				& Conver. rate  & --         & 2.42      &  2.10       &  2.14       \\
				\bottomrule
			\end{tabular}
		\end{center}
	\end{table}
	
	Next, we test the convergence rates in time. We fix the spatial mesh $(N_p,N_T)=(10370,20763)$ and run our Euler scheme and Crank-Nicolson scheme under different $\tau$. We then compute the convergence rate by the following formula at $t=0.24$
	\begin{equation}\label{rate_t}
		\mbox{Conver. rate}=\frac{\ln\left({\rm Error}(t;h_0,\tau)/{\rm Error}(t;h_0,\tau/2)\right)}{\ln 2}.
	\end{equation}
	\Cref{CT-MCF-t} presents ${\rm Error}(0.24;h_0,\tau)$ under various $\tau$ with $\omega=0.5$. It can be observed that the Euler scheme is first-order accurate, whereas the Crank-Nicolson scheme is second-order accurate. The results of  other $\omega$ are almost the same, so we omit them here.
	\begin{table}[htp!]
		\centering
		\footnotesize
		\caption{\Cref{ex1}: {${\rm Error}(0.24;h_0,\tau)$} under different $\tau$, $\tau_0=2\times 10^{-2}$, $\omega=0.5$.}\label{CT-MCF-t}
		\begin{tabular}{cccccc}
			\toprule
			& $\tau $		 &    $\tau_0/4$  & $\tau_0/8$     & $\tau_0/16$   & $\tau_0/32$ \\ 
			\cline{2-6}
			Euler            & Error 		& 7.7065e-02    & 4.5334e-02   & 2.4942e-02   & 1.2849e-02 \\ 
			& Conver. rate   & --           & 0.77          & 0.86          &  0.96  \\ 
			\toprule
			& $\tau $ 		 &  $\tau_0$    & $\tau_0/2$    & $\tau_0/4$    & $\tau_0/8$ \\
			\cline{2-6}
			Crank-Nicolson   & Error 	& 3.9705e-02    &1.4305e-02     & 3.8698e-03    & 2.9704e-04\\ 
			& Conver. rate   & --           & 1.47          & 1.89          &  3.70  \\ 
			\bottomrule
		\end{tabular}
	\end{table}

\end{example}

\begin{example}\label{ex2}
	Next, we consider the commonly used benchmark, where the initial surface is dumbbell-shaped given by the following local parametrization 
	\begin{equation}\label{dumbbell-06}
		\vx=\left[
		\begin{matrix}
			\cos\varphi\\
			(r\cos^2\varphi+1-r)\cos\theta\sin\varphi\\
			(r\cos^2\varphi+1-r)\sin\theta\sin\varphi
		\end{matrix}
		\right],  \qquad \theta\in[0,2\pi),\quad \varphi\in[0,\pi].
	\end{equation}
	
	The case with  $r=0.6$ is usually used to test the robustness of an algorithm in keeping mesh quality, e.g., in \cite{duan2024new,Elliott-Fritz-2017,hu2022evolving,duan2024energy,duan2025mesh}. The initial mesh we adopt here is completely the same as that in \cite{duan2024energy,duan2025mesh}, which  contains $N_T=4192$ plane triangles with $N_p=2098$ vertices, see \Cref{fig-dumbbell-06}(a).  \Cref{fig-dumbbell-06}(b)--(c) present the  numerical surfaces at $t=0.06$ and $t=0.088$ obtained by our Crank-Nicolson scheme with $\omega=1$ and $\tau=1\times 10^{-4}$. To capture the evolution near the singularity, we reduce the time step size to $\tau=1\times 10^{-6}$ for $t\ge 0.0879$. One can see the mesh is preserved as well as the initial mesh  during the evolution. \Cref{fig-dumbbell-06-Energy}(a) shows $|\Gamma_h^m|$  with respect to time which demonstrates the energy-stable property. \Cref{fig-dumbbell-06-Energy}(b) presents $\sigma_{\max}$ for our Crank-Nicolson scheme. We also run our code under $\omega=0.5$ and $0.1$ and get similar results, so we omit them here.  
	\begin{figure}[!htp]
		\centering
		\subfigure[Initial surface, $|\Gamma_h|=6.123$]{\includegraphics[width=.32\textwidth,trim=20 200 20 230, clip]{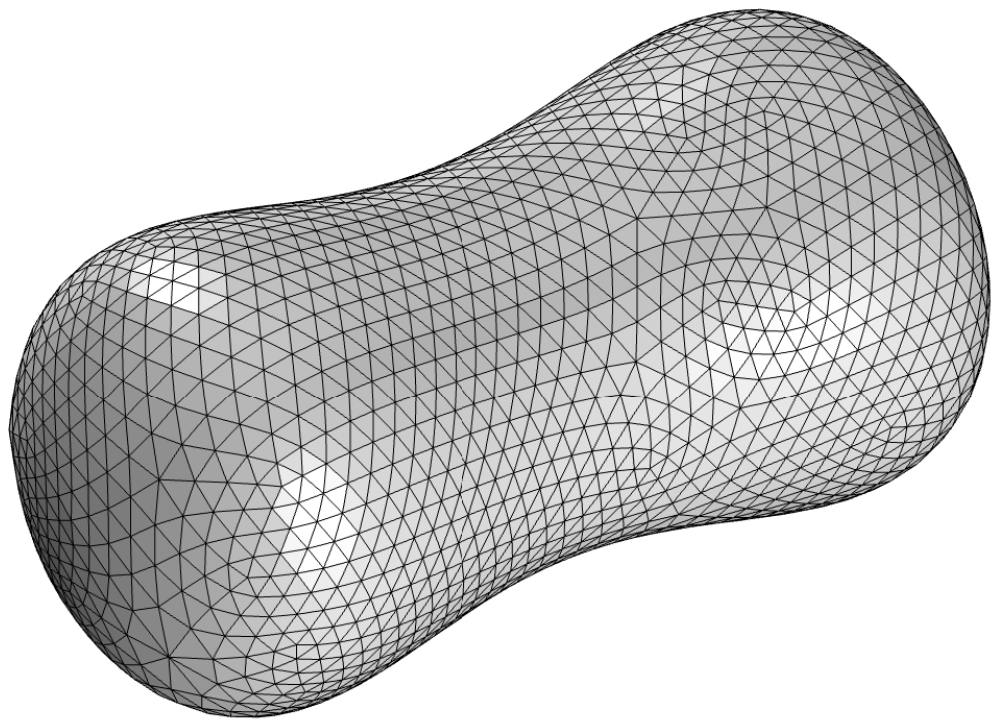}}
		\subfigure[$t=0.06$, $|\Gamma_h|=1.922$]{\includegraphics[width=.32\textwidth,trim=20 200 20 250, clip]{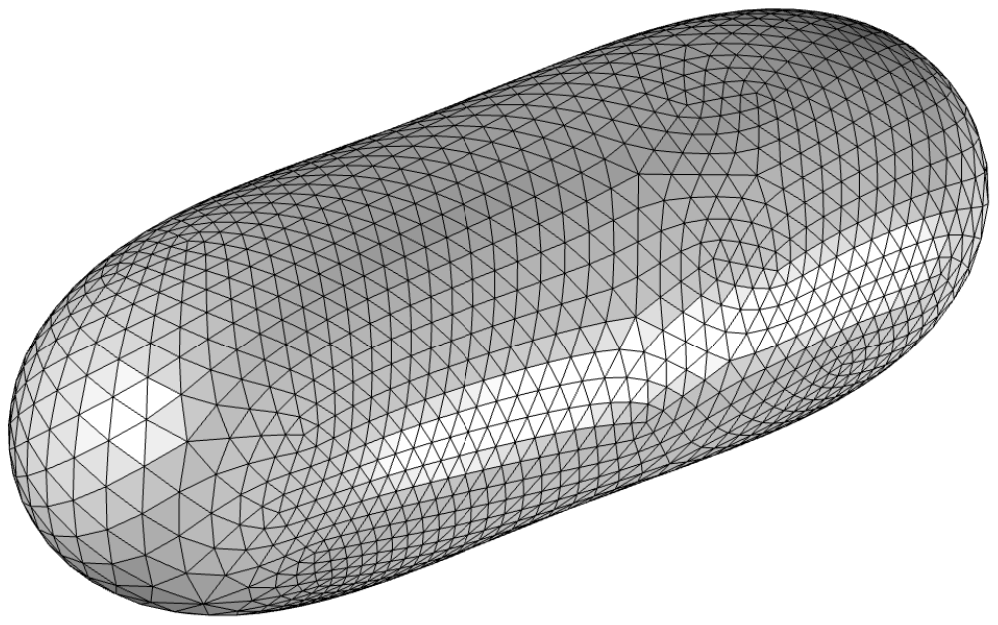}}
		\subfigure[$t=0.088$, $|\Gamma_h|=2.436\times 10^{-5}$]{\includegraphics[width=.32\textwidth,trim=0 160 0 190, clip]{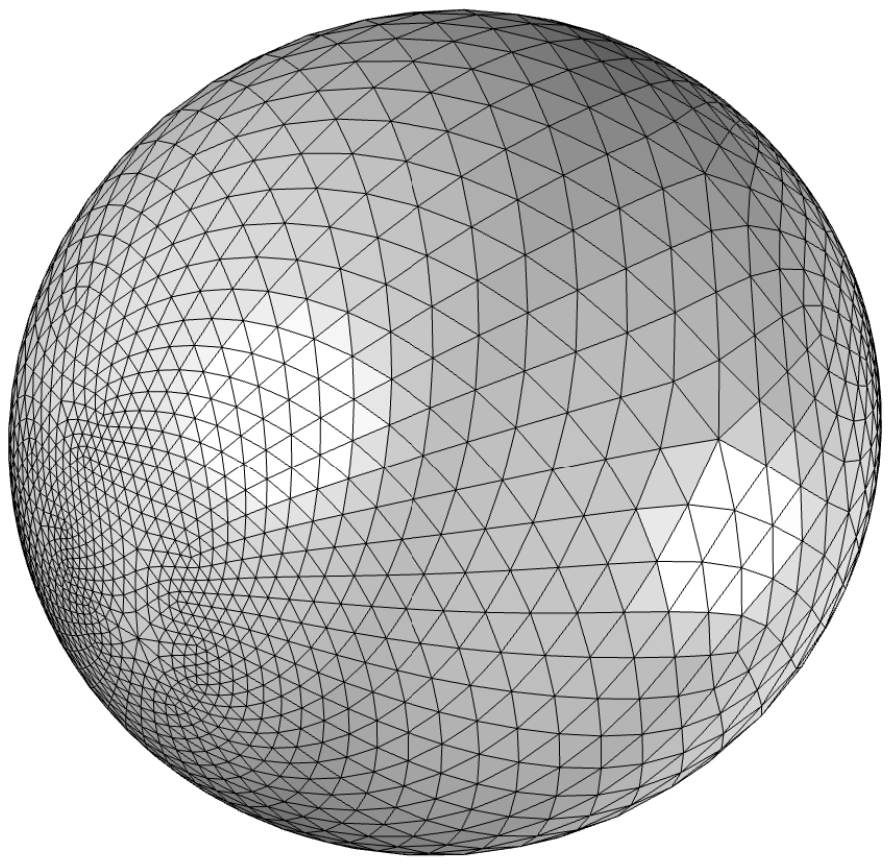}}
		\caption{Numerical surfaces at different times for \Cref{ex2}, the images are rescaled.}
		\label{fig-dumbbell-06}
	\end{figure}
	\begin{figure}[!htp]
		\centering
		\subfigure[Evolution of surface areas]{\includegraphics[width=.42\textwidth,trim=0 200 0 200, clip]{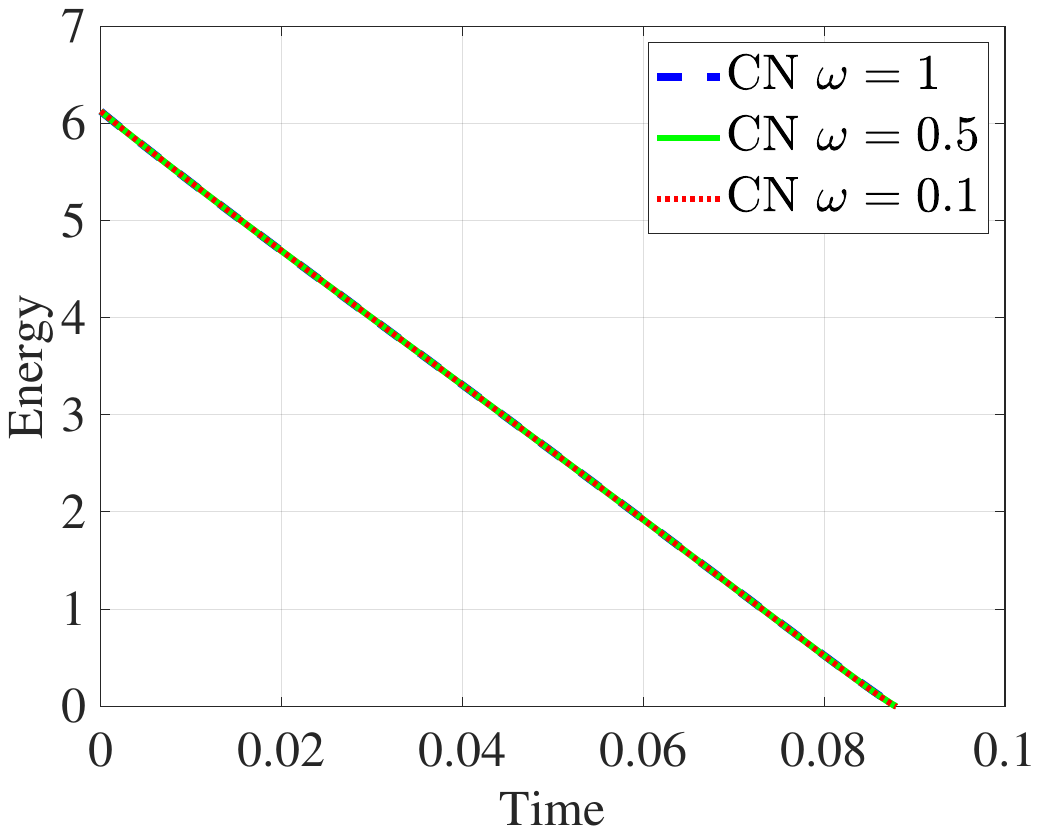}}
		\subfigure[The mesh quality index $\sigma_{\max}$]{\includegraphics[width=.42\textwidth,trim=0 200 0 200, clip]{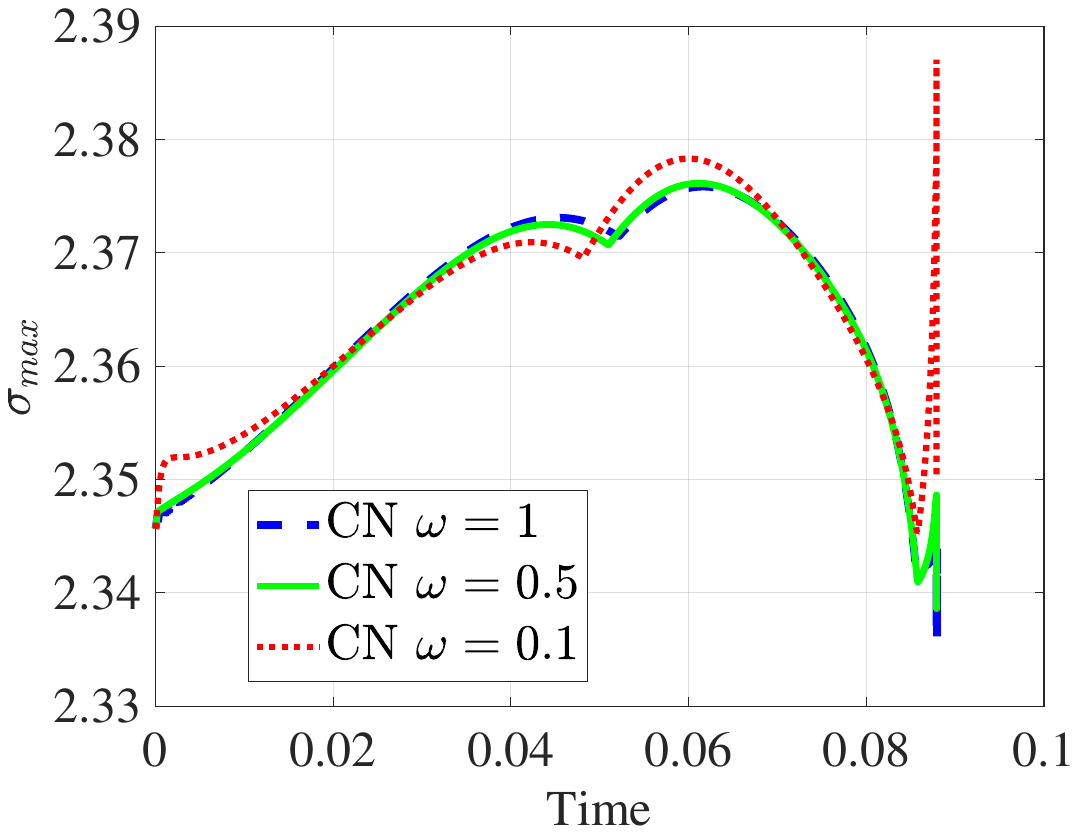}}
		\caption{Surface areas and $\sigma_{\max}$ with respect to time  for \Cref{ex2}, $r=0.6$.}
		\label{fig-dumbbell-06-Energy}
	\end{figure}
	
	We next change $r$ from $0.6$ to $0.7$ to show the effect of $\omega$ when singularities develop during the evolution. We adopt the same initial mesh as in \cite[Example 5.2]{duan2024energy}, which has $N_T=4352$ plane triangles and $N_p=2178$ vertices. We run our Crank-Nicolson scheme with $\tau=1\times 10^{-4}$ for different values of $\omega$ and show the numerical surfaces at different times in \Cref{fig-dumbbell-07}. One can observe that for $\omega=0.5$ the tangential motion drives the vertices to move toward the neck so as to maintain conformality. As a result, the mesh becomes coarse at the two ends, which is similar to the results obtained using the DeTurck method in \cite[Fig. 14]{Elliott-Fritz-2017}. We plot $\sigma_{\max}$ for $\omega=0,0.1,0.5,1$ in \Cref{fig-dumbbell-07-Energy}(b), from which we observe that the mesh quality index $\sigma_{\max}$ does not vary significantly among these values of $\omega$ before the singularity developed. In fact, $\sigma_{\max}\approx 5.84,6.40,7.14,8.09$ correspond to $\omega=1,0.5,0.1,0$, respectively, at $t=0.0565$. \Cref{fig-dumbbell-07-Energy}(a) presents the surface area for different values of $\omega$, which again demonstrates the energy stability. Moreover, we observe that a smaller $\omega$ leads to faster evolution as the singularity is approached.
	\begin{figure}[!htp]
		\centering
		\subfigure[Initial surface]{\includegraphics[width=.32\textwidth,trim=0 200 0 200, clip]{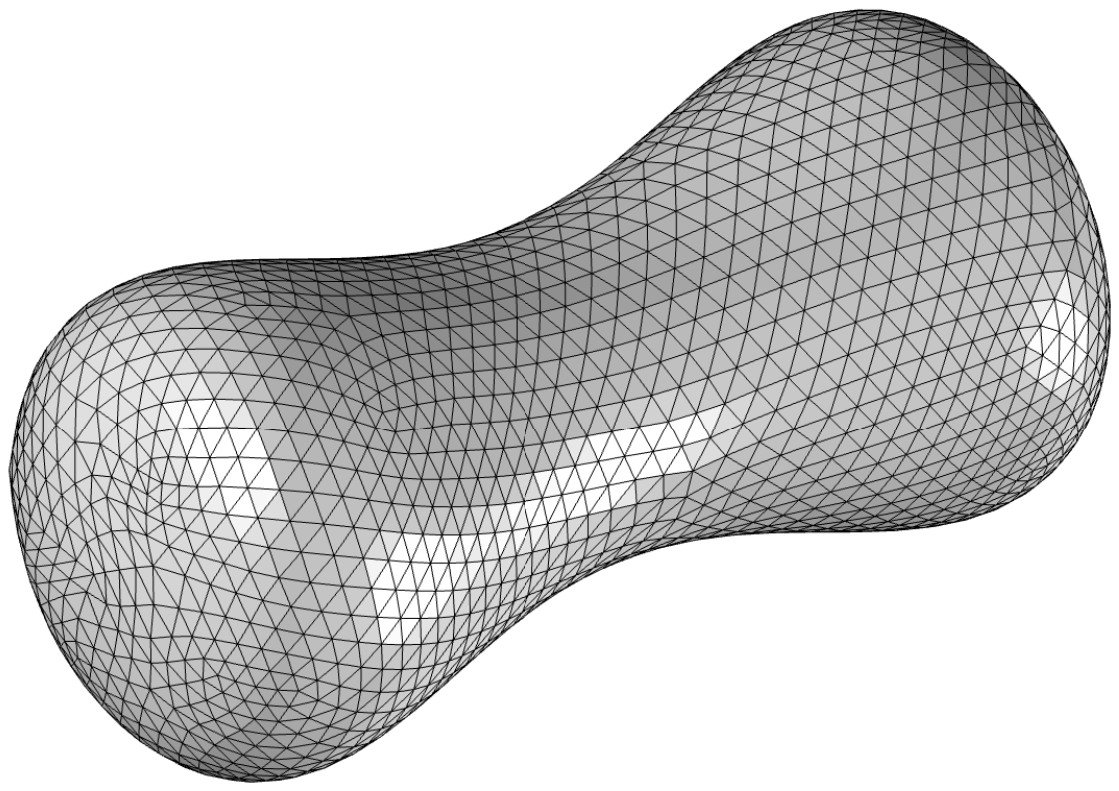}}
		\subfigure[$t=0.03$, $\omega=0 $]{\includegraphics[width=.32\textwidth,trim=0 200 0 240, clip]{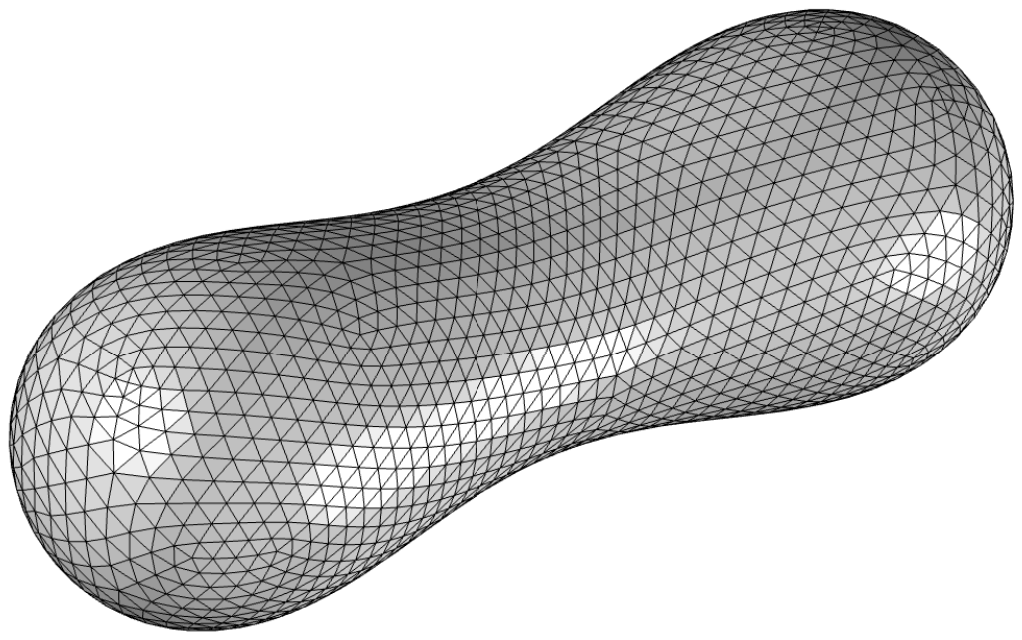}}
		\subfigure[$t=0.05$, $\omega=0 $]{\includegraphics[width=.32\textwidth,trim=0 180 0 180, clip]{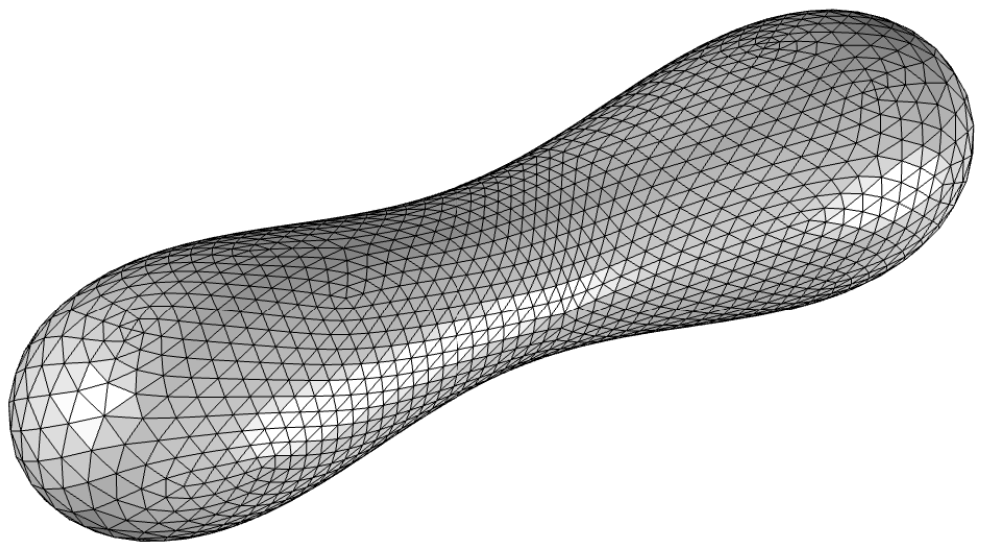}}
		\subfigure[$t=0.057$, $\omega=0 $]{\includegraphics[width=.32\textwidth,trim=30 240 30 240, clip]{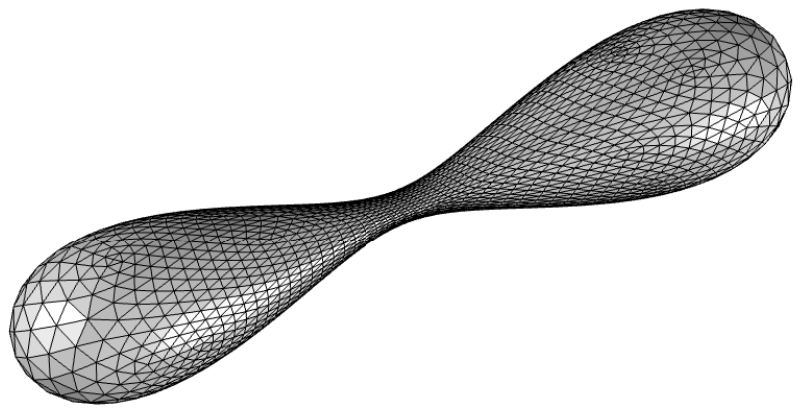}}
		\subfigure[$t=0.059$, $\omega=0 $]{\includegraphics[width=.32\textwidth,trim=30 240 30 230, clip]{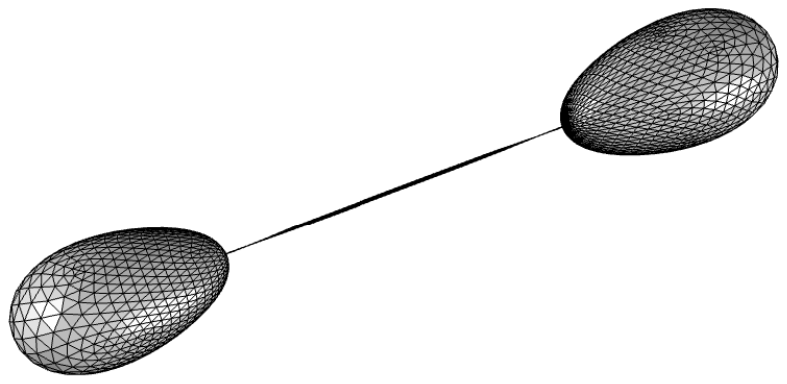}}
		\subfigure[$t=0.06$, $\omega=0 $]{\includegraphics[width=.32\textwidth,trim=30 240 30 230, clip]{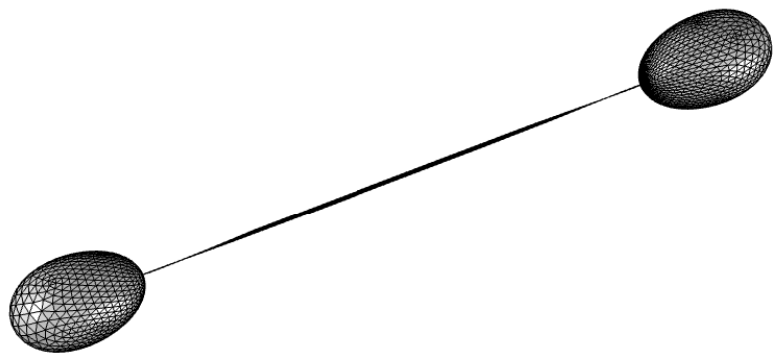}}
		\subfigure[$t=0.057$, $\omega=0.5 $]{\includegraphics[width=.32\textwidth,trim=0 240 0 220, clip]{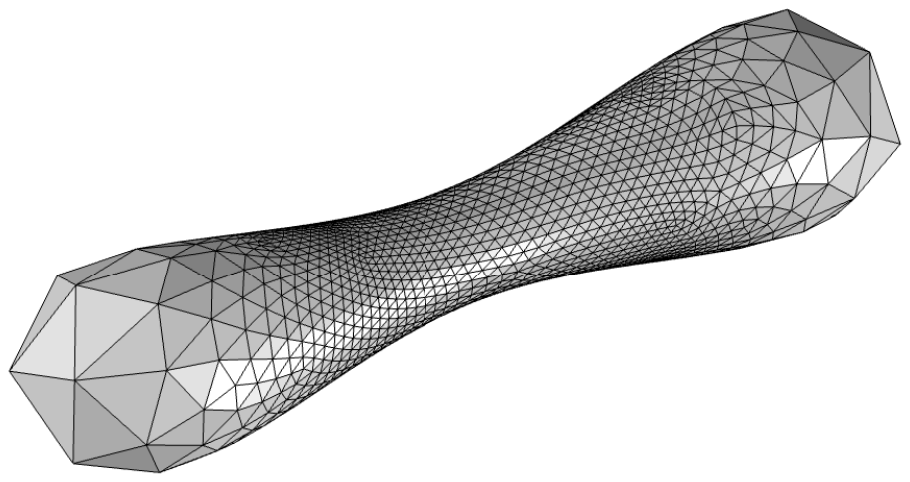}}
		\subfigure[$t=0.059$, $\omega=0.5 $]{\includegraphics[width=.33\textwidth,trim=0 240 0 220, clip]{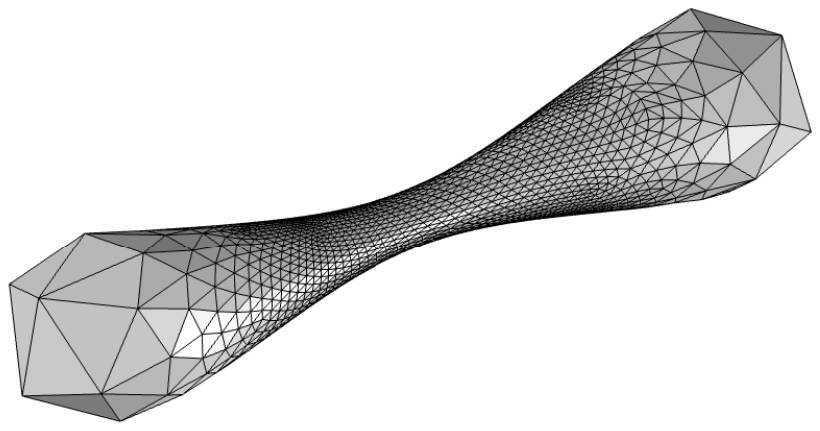}}
		\subfigure[$t=0.06$, $\omega=0.5 $]{\includegraphics[width=.28\textwidth,trim=0 240 0 220, clip]{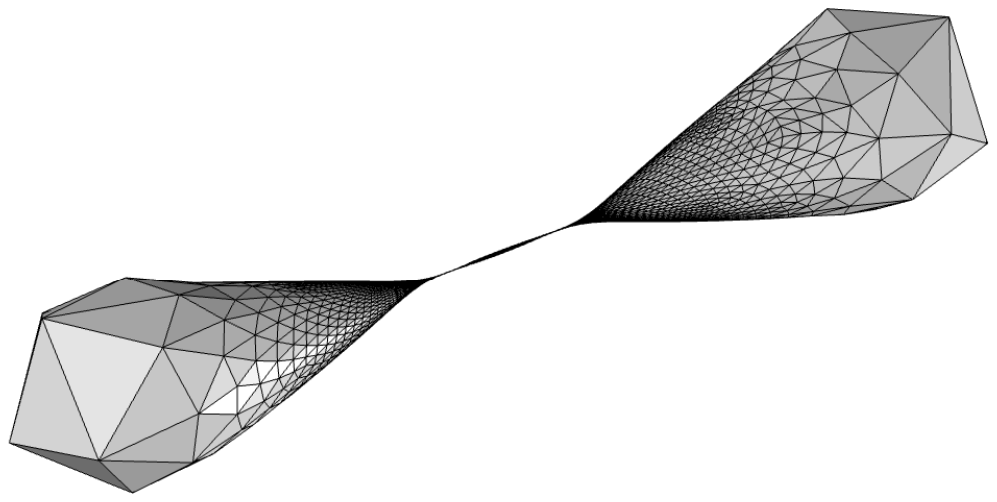}}
		\caption{Numerical surfaces at different times for \Cref{ex2}, the images are rescaled.}
		\label{fig-dumbbell-07}
	\end{figure}
	
	\begin{figure}[!htp]
		\centering
		\subfigure[Evolution of surface areas]{\includegraphics[width=.42\textwidth,trim=0 200 0 200, clip]{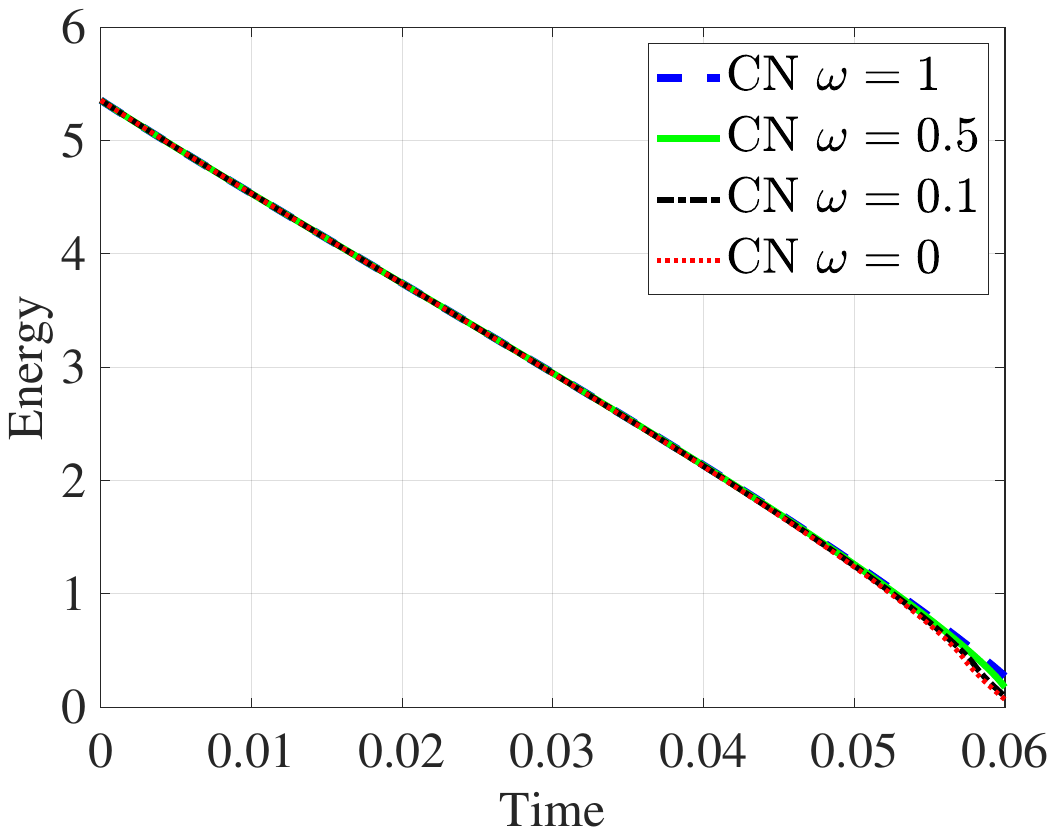}}
		\subfigure[The mesh quality $\sigma_{\max}$]{\includegraphics[width=.42\textwidth,trim=0 200 0 200, clip]{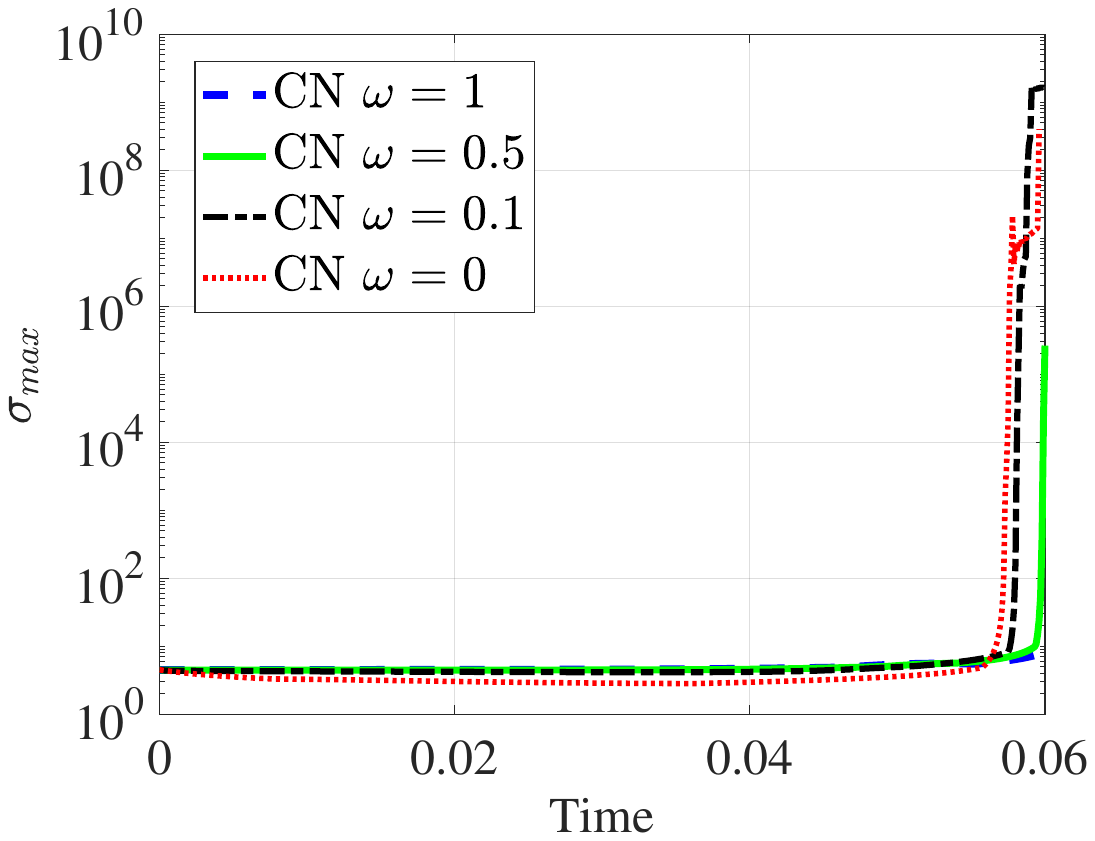}}
		\caption{Surface areas and $\sigma_{\max}$ with respect to time  for \Cref{ex2}, $r=0.7$.}
		\label{fig-dumbbell-07-Energy}
	\end{figure}
\end{example}

\subsection{Surface diffusion}
\begin{example}\label{ex3}
	Note that, in contrast to \eqref{eqn-Euler} and \eqref{eqn-CN}, the formulations \eqref{Euler-SD} and \eqref{CN-SD} involve boundary terms; therefore, we test the convergence rates only for the schemes \eqref{Euler-SD} and \eqref{CN-SD}. The initial surface is half of an ellipsoid given by
	\[
	\frac{x^2}{2^2}+\frac{y^2}{1.5^2}+\frac{z^2}{1^2}=1,\qquad z\ge 0.
	\]
	We set $\alpha=\frac{2}{3}\pi$ and $\beta=100$. To test the spatial convergence rate, we construct a family of meshes for $\Gamma_0$, denoted by $\{\mathcal{T}_l\}_{l=0}^4$. The mesh $\mathcal{T}_l$ is obtained by refining each triangle of $\mathcal{T}_{l-1}$ into four smaller triangles and then lifting the new nodes onto the surface. The initial mesh $\mathcal{T}_0$ contains $205$ triangles and $116$ vertices. The numbers of vertices and triangles for $\{\mathcal{T}_l\}_{l=1}^4$ are listed below:
	\begin{equation*}
		(N_p,N_T)=(436,820),\quad (1691,3280),\quad (6661,13120),\quad (26441,52480).
	\end{equation*}
	Let $\Gamma_h^m(\mathcal{T}_l;\tau)$ denote the numerical surface obtained by either the Euler or the Crank-Nicolson scheme when $\mathcal{T}_l$ is chosen as the initial mesh and $\tau$ is the time stepsize. We define the spatial error as
	\begin{equation}\label{eqn-err-def}
		{\rm Error}_s(t_m;l,\tau)=\mathcal{H}[\Gamma_h^m(\mathcal{T}_{l};\tau),\Gamma_h^m(\mathcal{T}_{l+1};\tau)].
	\end{equation}
	The convergence rate is then computed using \eqref{rate_s}. \Cref{CT-SD-h} presents ${\rm Error}_s(0.2;l,\tau)$ for $l=0,1,2,3$, where the numerical surfaces are obtained using our Crank-Nicolson scheme with $\tau=4\times 10^{-4}$. One can observe second-order convergence in \Cref{CT-SD-h}.
	
	\begin{table}[htp!]
		\centering
		\footnotesize
		\caption{\Cref{ex3}: ${\rm Error}_s(0.2;l,5\times 10^{-4})$ under different $\omega$ for $l=0,1,2,3$.}\label{CT-SD-h}
		\begin{center}
			\begin{tabular}{cccccc}
				\toprule
				$\omega$            & $h_0$   	& 0.5949     &   0.3006     &  0.1507     &   0.0754	        \\
				\toprule
				\multirow{2}{*}{1}  & Error     & 2.9624e-02  & 9.6222e-03  & 2.6107e-03  & 6.5810e-04       \\
				& Conver. rate      & --          & 1.65                  &  1.89              &  1.99      \\
				\toprule
				\multirow{2}{*}{0.5}& Error     & 2.7709e-02  & 9.3062e-03  & 2.4633e-03  &  6.1727e-04  \\
				& Conver. rate      & --            & 1.60                 &  1.93        	     &  2.00      \\
				\toprule
				\multirow{2}{*}{0}    & Error & 3.2311e-02    & 6.9279e-03  & 1.7656e-03  & 5.3501e-04    \\
				& Conver. rate      & --        & 2.26       &  1.98       &  1.73      \\
				\bottomrule
			\end{tabular}
		\end{center}
	\end{table}
	
	To test the temporal convergence rate, we fix the initial mesh and run our schemes with different time stepsizes. We compute the error between $\Gamma_h^{m}(\mathcal{T}_4;\tau)$ and $\Gamma_h^{m}(\mathcal{T}_4;\tau')$ by
	\begin{equation}\label{eqn-err-t-def}
		{\rm Error}_t(t_m;4,\tau)
		= \mathcal{H}\big[\Gamma_h^m(\mathcal{T}_4;\tau),\Gamma_h^m(\mathcal{T}_4;\tau')\big].
	\end{equation}
	Note that the error ${\rm Error}_t(t_m;4,\tau)$ defined above also involves spatial errors. For the Euler scheme, we set $\tau'=\frac{1}{2}\tau$. For the Crank-Nicolson scheme, since the error decays rapidly when the time stepsize is reduced, we set $\tau'=\frac{k}{k+1}\tau$ with $k\in\mathbb{N}^+$. \Cref{CT-SD-t-Euler} and \Cref{CT-SD-t-CN} show the results, where  $\mathcal{T}_4$ is chosen as the initial triangulation. For the Crank-Nicolson scheme, we present results for $\omega=1,0.5,0$, while for the Euler scheme we show results for $\omega=1,0.5,0.1$, since our Euler scheme reduces to the BGN scheme when $\omega=0$. The convergence rate is computed using \eqref{rate_t}.  We observe that our Euler scheme is of order $\mathcal{O}(\tau)$, while our Crank-Nicolson scheme converges at the rate about $\mathcal{O}(\tau^2)$. 
	\begin{table}[htp!]
		\centering
		\footnotesize
		\caption{\Cref{ex1}: ${\rm Error}_t(0.2;4,\tau)$ of our Euler scheme, $\tau_0=2.5\times 10^{-2}$.}\label{CT-SD-t-Euler}
		\begin{tabular}{cccccc}
			\toprule
			$\omega$             & $\tau $	&  $\tau_0/8$    & $\tau_0/16$    & $\tau_0/32$   & $\tau_0/64$ \\
			\toprule
			\multirow{2}{*}{1}   & Error   & 6.2942e-03      & 3.4403e-03     & 1.8326e-03    & 9.5488e-04 \\ 
			& Conver. rate    	 & --      & 0.87            & 0.91           &  0.94 \\ 
			\cline{1-6}
			\multirow{2}{*}{0.5} & Error   & 5.9109e-03      & 3.2586e-03    & 1.7569e-03     & 9.2820e-04 \\ 
			& Conver. rate    	 & --      & 0.86            & 0.89          &  0.92  \\ 
			\cline{1-6}
			\multirow{2}{*}{0.1} & Error   & 4.1325e-03     & 2.6379e-03     & 1.4951e-03     & 8.1868e-04\\ 
			& Conver. rate    	 & --      & 0.65          & 0.82         &  0.87  \\ 		
			\bottomrule
			
		\end{tabular}
	\end{table}
	%
	\begin{table}[htp!]
		\centering
		\footnotesize
		\caption{\Cref{ex1}: ${\rm Error}_t(0.2;4,\tau)$ of our Crank-Nicolson scheme, $\tau_0=2.5\times 10^{-2}$.}\label{CT-SD-t-CN}
		\begin{tabular}{cccccc}
			\toprule
			$\omega$             & $\tau $	&  $\tau_0/3$    & $\tau_0/4$    & $\tau_0/5$   & $\tau_0/6$       \\
			\toprule
			\multirow{2}{*}{1}   & Error 	    & 2.1842e-03 &  1.3876e-03  & 9.6861e-04  & 7.1915e-04   \\ 
			& Conver. rate    & --          & 1.57          & 1.61        &  1.63 \\ 
			\cline{1-6}
			\multirow{2}{*}{0.5}   & Error 	    &  1.7920e-03  & 1.1416e-03  & 8.0030e-04 &  5.9631e-04 \\ 
			& Conver. rate    & --          & 1.57         & 1.59        &  1.61\\ 
			\cline{1-6}
			\multirow{2}{*}{0}   & Error 	    & 7.0215e-04  & 4.0194e-04   &2.6267e-04   &1.8956e-04\\ 
			& Conver. rate    & --          & 1.93          & 1.91        &  1.79 \\ 
			\bottomrule
		\end{tabular}
	\end{table}
	
	We adopt $\mathcal{T}_2$ as the initial mesh and run our Crank-Nicolson scheme with $\tau=5\times 10^{-3}$ under different $\omega$.  \Cref{fig-half-ellipsoid}(a)--(c) present the numerical surfaces under $\omega=1$ at $t=0,0.5$ and $1$, respectively.  Denote by 
	\begin{equation}\label{Energy-SD}
		E_h^{m}:=|\Gamma_h^m| -\cos\alpha|\Omega_h^{m}|
	\end{equation}
	the energy at $t=t_m$, then for $\omega=0.1,0.5$ and $1$, $\max_m(E_h^m-E_h^{m-1})\approx-3.31\times 10^{-5}$, $-3.54\times 10^{-5}$ and $-3.79\times 10^{-5}$, respectively. And in all the cases, $\max_m{\sigma_{\max}}\approx 3.1$. In fact, the initial mesh yields $\sigma_{\max}\approx 3.0$. 
	
	\begin{figure}[!htp]
		\centering
		\subfigure[Initial surface]{\includegraphics[width=.31\textwidth,trim=0 200 0 200, clip]{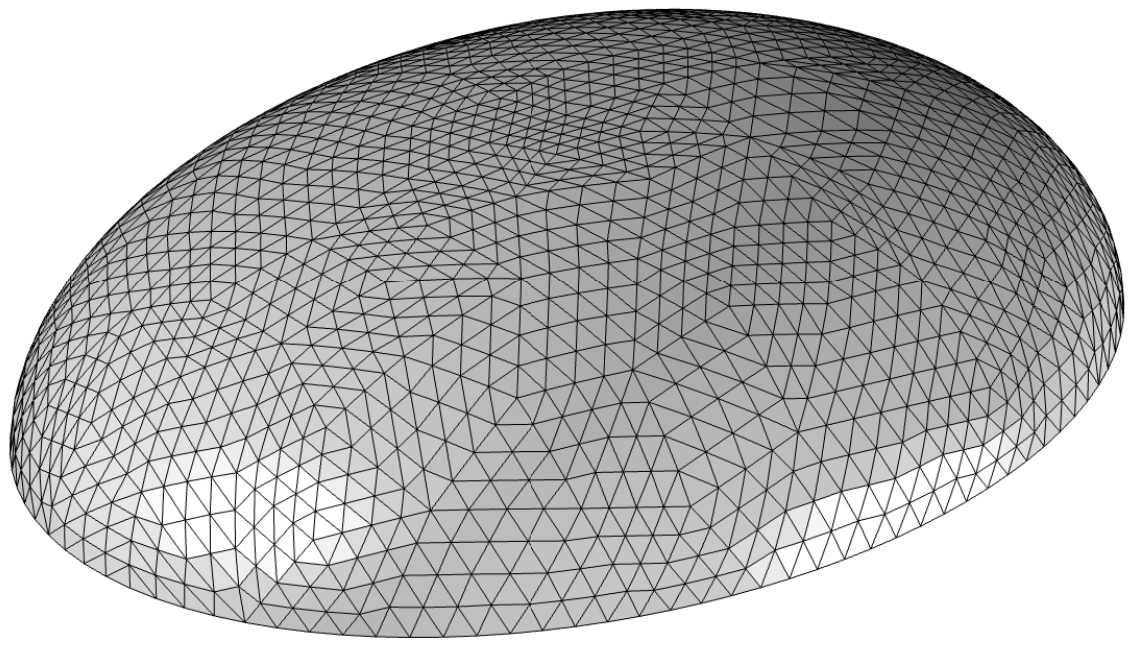}}
		\subfigure[$t=0.5$]{\includegraphics[width=.31\textwidth,trim=0 180 0 180, clip]{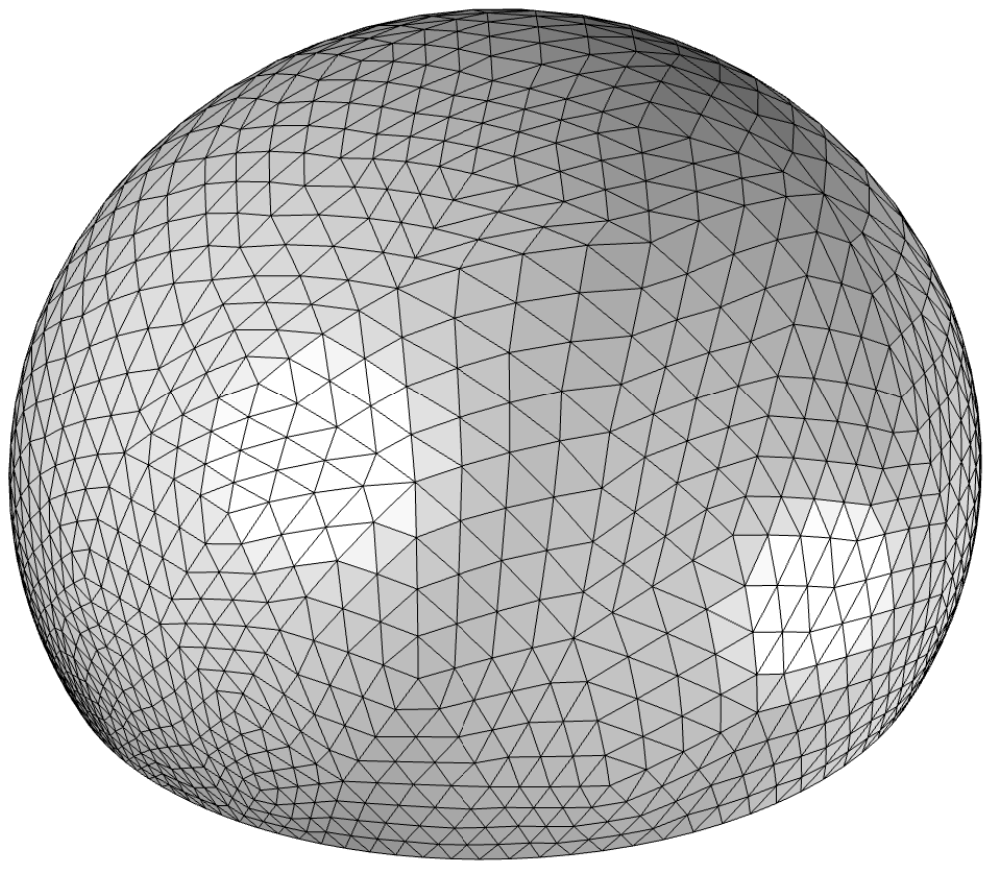}}
		\subfigure[$t=1$]{\includegraphics[width=.31\textwidth,trim=0 160 0 180, clip]{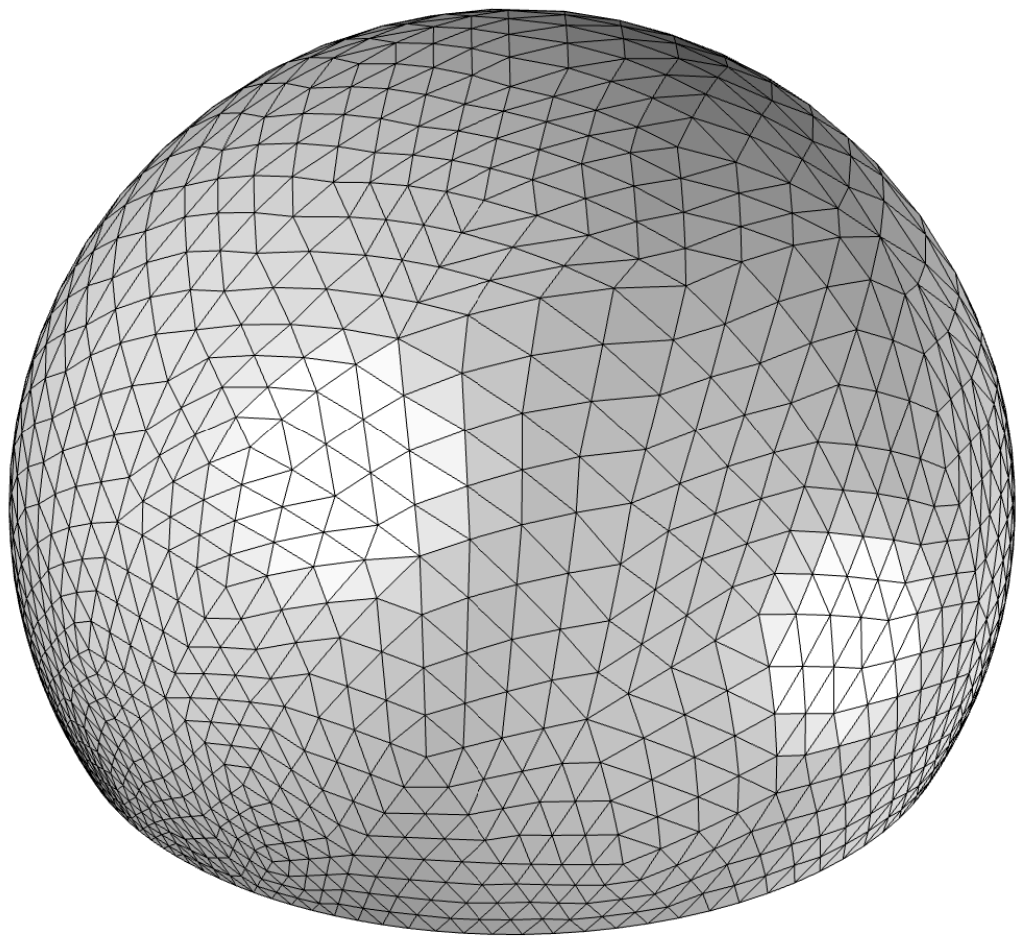}}
		\caption{Numerical surfaces at different times for \Cref{ex3}, the images are rescaled.}
		\label{fig-half-ellipsoid}
	\end{figure}
\end{example}

\begin{example}\label{ex-half-box}
	In this example, we consider the solid-state dewetting model with $\alpha=\frac{2}{3}\pi$ and $\beta=100$. The initial surface is a $3.2\times 3.2\times 0.1$ cuboid island. This example was also considered in \cite[Example 5.1]{bao2023}. The initial mesh contains $N_T=7604$ triangles and $N_p=3911$ vertices, with $216$ of them on the boundary. Since the initial surface has sharp corners, we use a graded mesh $t_m=T(m/200)^2$ with $T=0.1$. The numerical solutions at the first $4$ steps are obtained from our Euler scheme; otherwise, the mesh quality index $\sigma_{\max}$ would oscillate (although the energy still decays smoothly). Of course, one can use a finer stepsize to suppress this phenomenon. \Cref{fig-half-box}(b)--(f) show the numerical surfaces at different times. The results are consistent with those given in \cite[Fig. 5.3]{bao2023}: the surface shrinks to a cross-shaped geometry since the four edges move faster than the four corners at the beginning, and after some time the corners dominate the evolution. Finally, it evolves into a spherical shape as equilibrium. The mesh quality with respect to time is plotted in \Cref{fig-half-box-Energy}(b), which again demonstrates the ability of our method to preserve mesh quality. \Cref{fig-half-box-Energy}(a) shows the energy given in \eqref{Energy-SD}, which verifies the energy stability.
	\begin{figure}[!htp]
		\centering
		\subfigure[Initial surface]{\includegraphics[width=.32\textwidth,trim=0 200 0 200, clip]{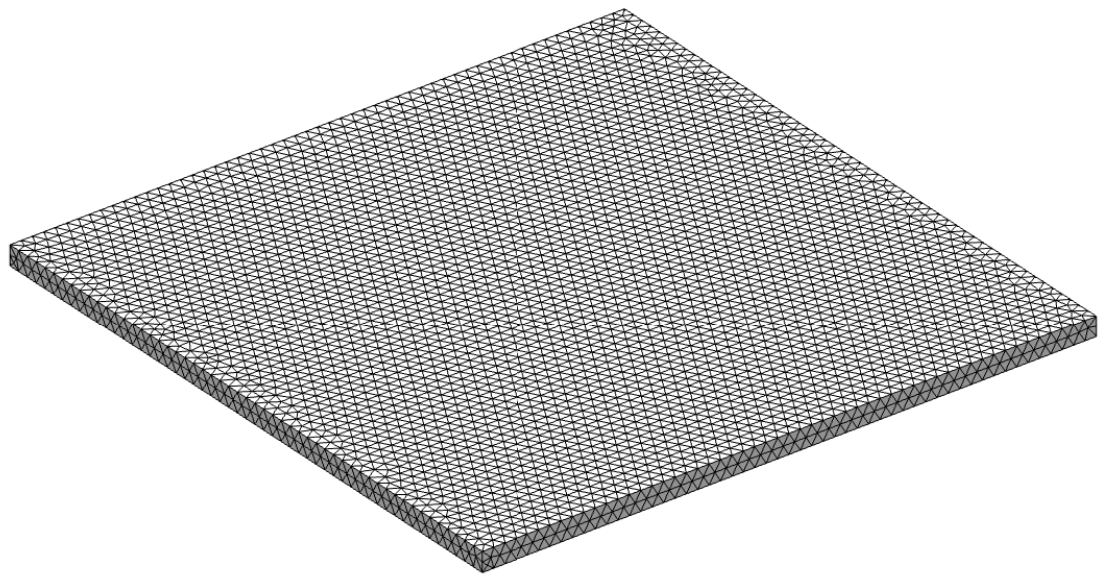}}
		\subfigure[$t=0.004$]{\includegraphics[width=.32\textwidth,trim=0 200 0 200, clip]{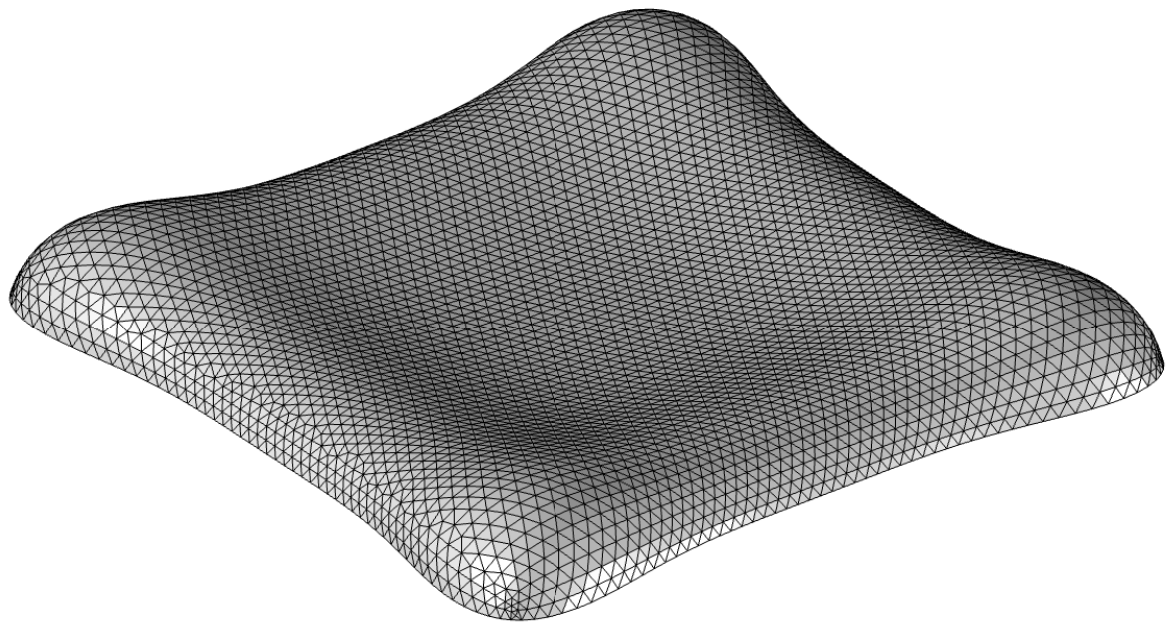}}
		\subfigure[$t=0.009$]{\includegraphics[width=.32\textwidth,trim=0 160 0 140, clip]{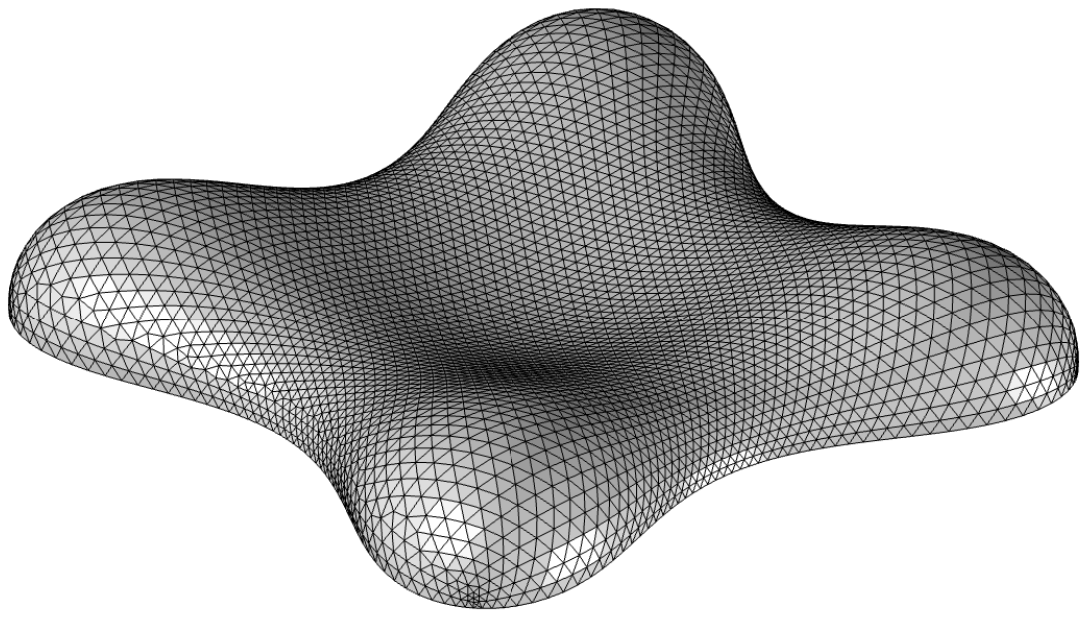}}
		\subfigure[$t=0.02025$]{\includegraphics[width=.32\textwidth,trim=0 160 0 140, clip]{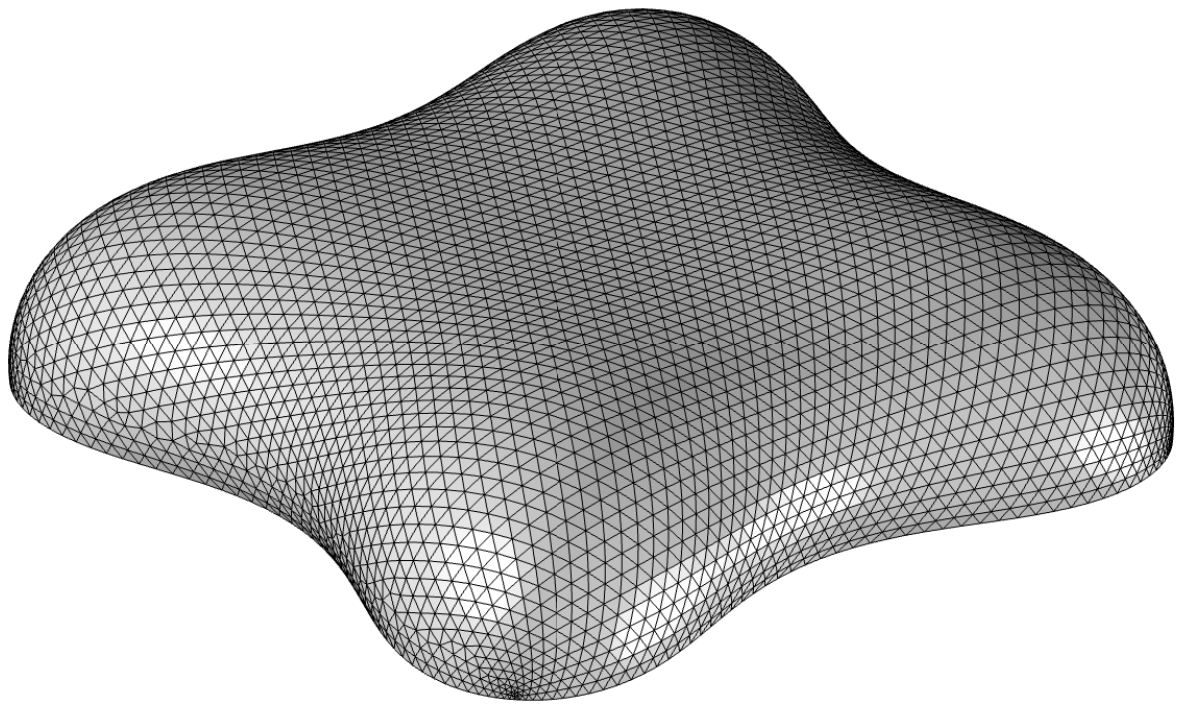}}
		\subfigure[$t=0.036$]{\includegraphics[width=.32\textwidth,trim=0 160 0 140, clip]{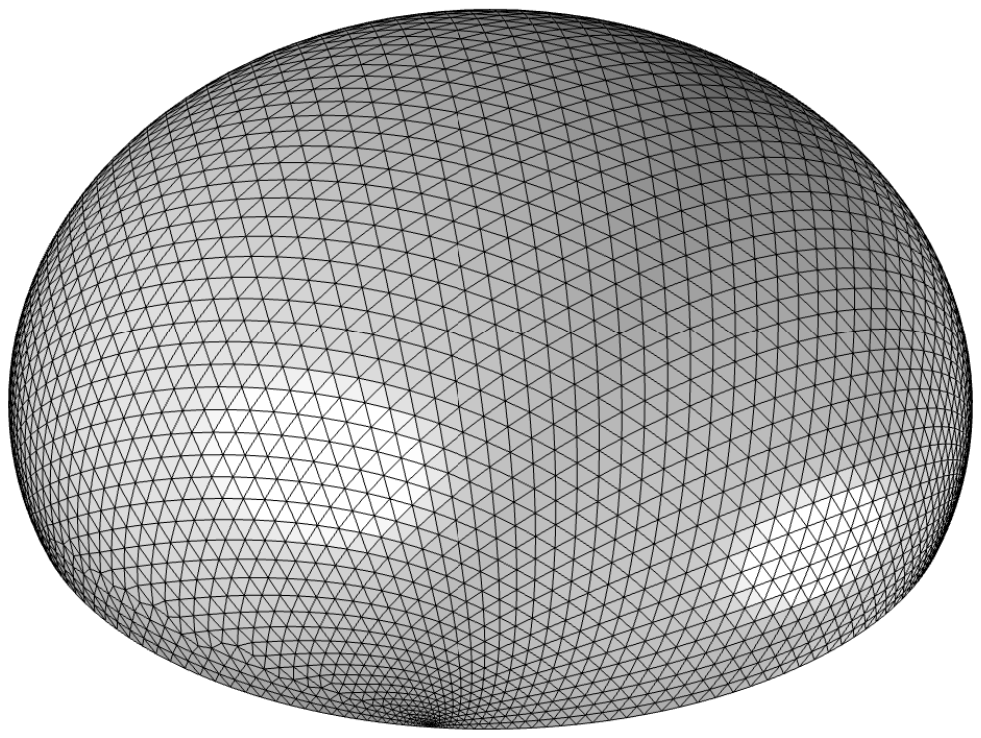}}
		\subfigure[$t=0.1$]{\includegraphics[width=.32\textwidth,trim=0 160 0 140, clip]{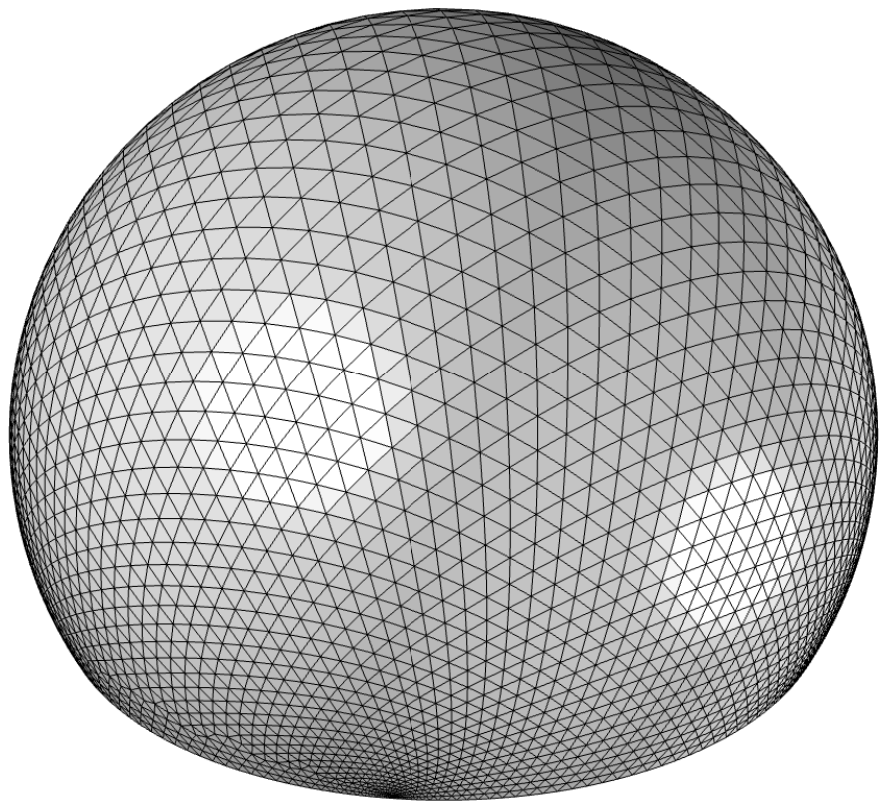}}
		\caption{Numerical surfaces at different times for \Cref{ex-half-box}, the images are rescaled.}
		\label{fig-half-box}
	\end{figure}
	
	\begin{figure}[!htp]
		\centering
		\subfigure[Evolution of surface areas]{\includegraphics[width=.42\textwidth,trim=0 200 0 200, clip]{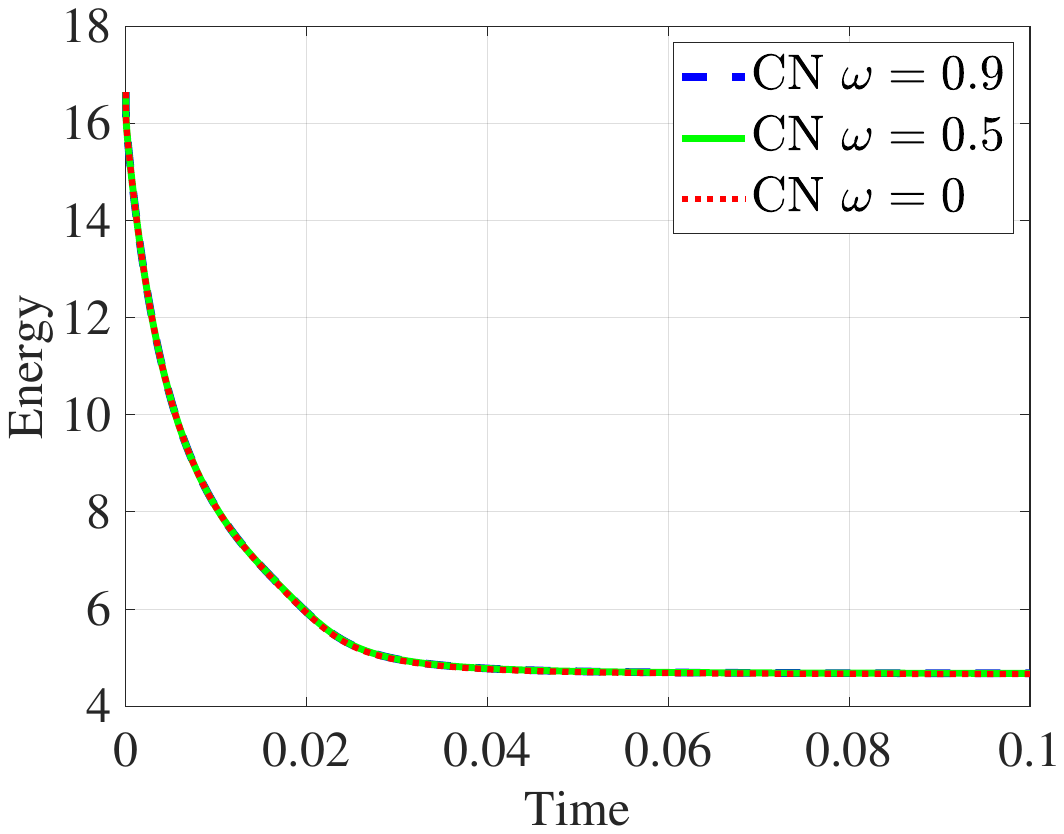}}
		\subfigure[The mesh quality index $\sigma_{\max}$]{\includegraphics[width=.42\textwidth,trim=0 200 0 200, clip]{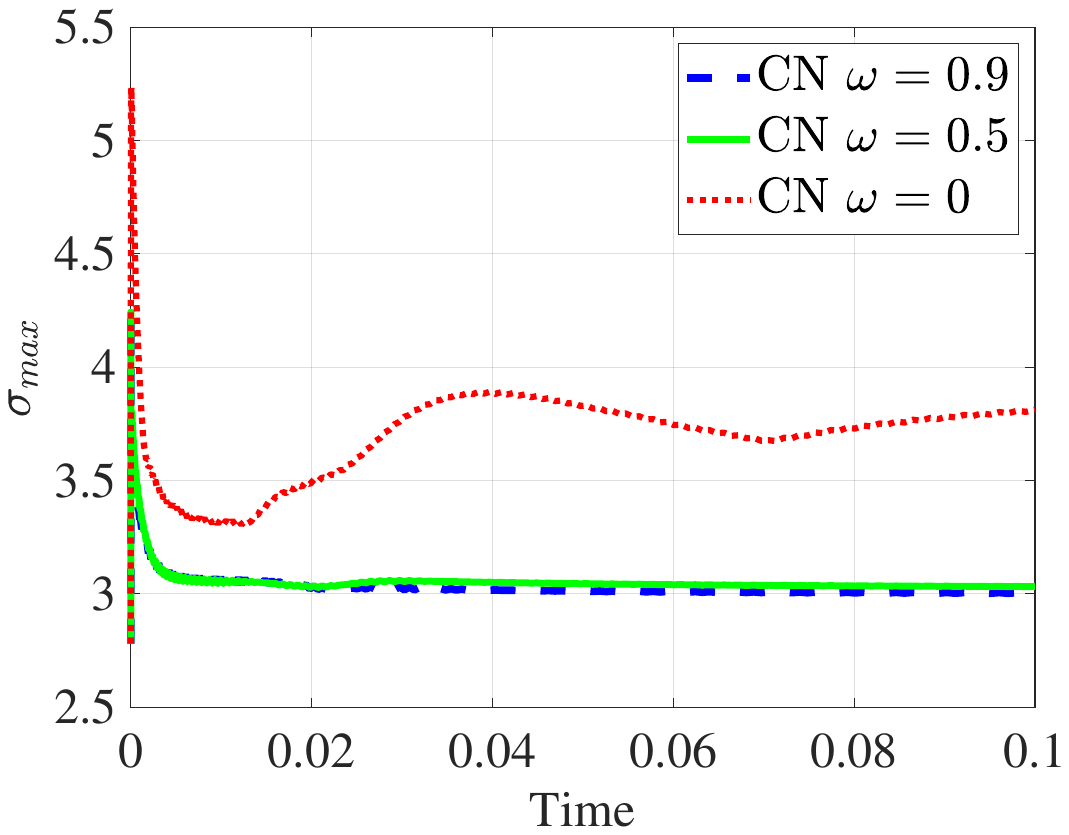}}
		\caption{Surface areas and $\sigma_{\max}$ with respect to time  for \Cref{ex-half-box}.}
		\label{fig-half-box-Energy}
	\end{figure}
	
\end{example}

\begin{figure}[!htp]
	\centering
	\subfigure[Initial surface]{\includegraphics[width=.32\textwidth,trim=0 250 0 200, clip]{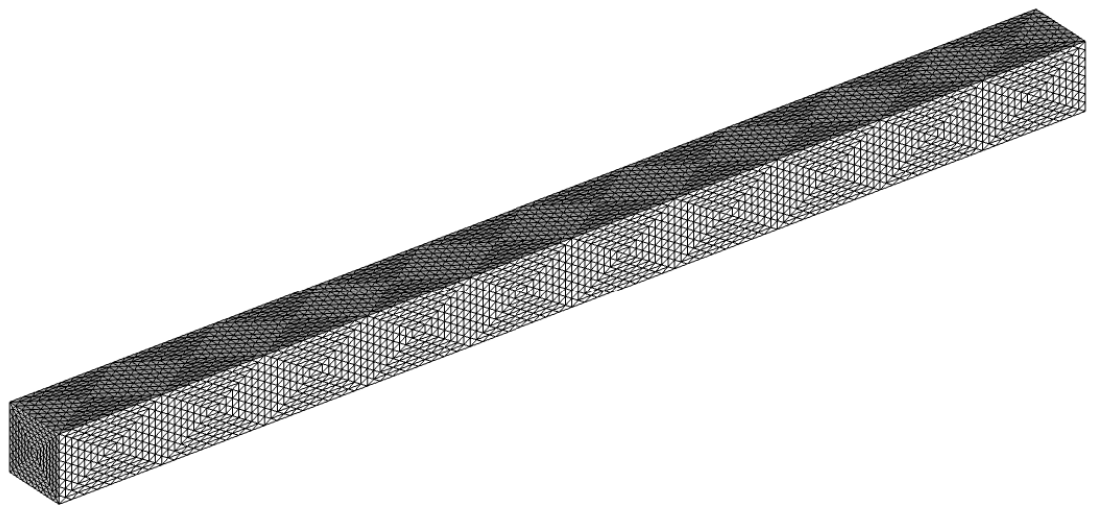}}
	\subfigure[Our CN, $t\approx 0.4612$]{\includegraphics[width=.32\textwidth,trim=0 250 0 250, clip]{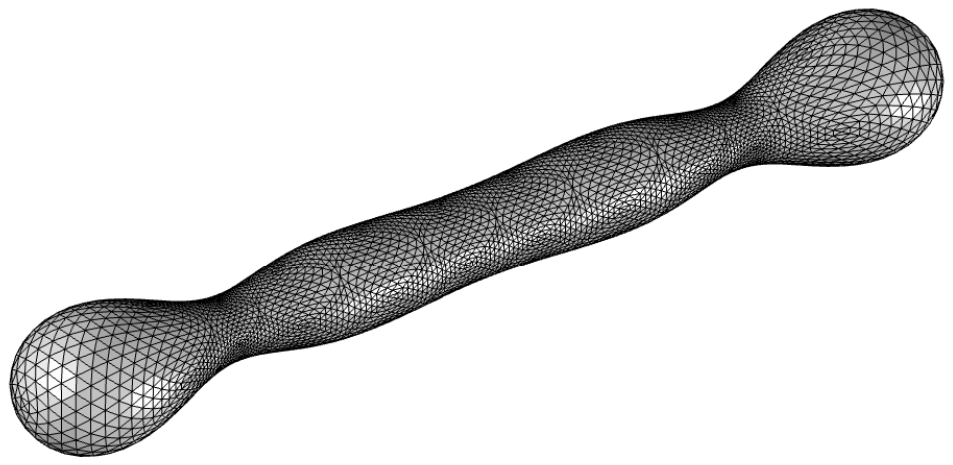}}
	\subfigure[Our CN, $t\approx 0.6367$]{\includegraphics[width=.32\textwidth,trim=0 250 0 250, clip]{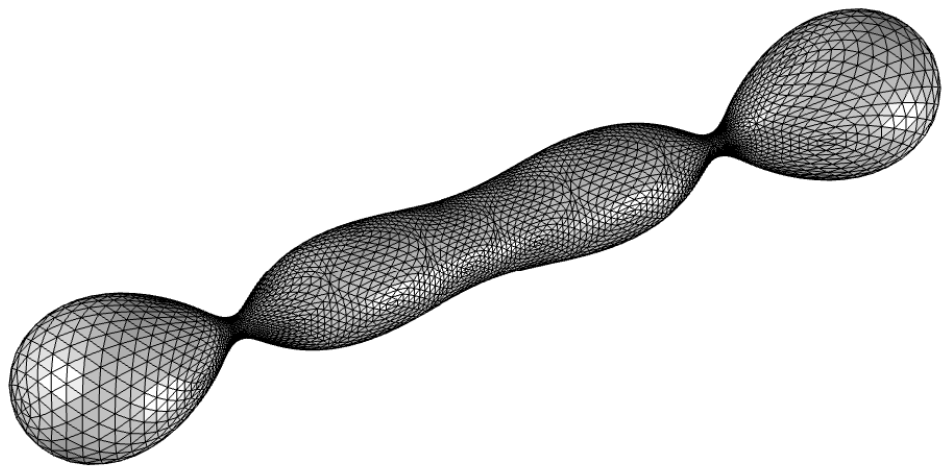}}
	
	\subfigure[Our CN, $t\approx 0.6720$]{\includegraphics[width=.32\textwidth,trim=60 270 60 250, clip]{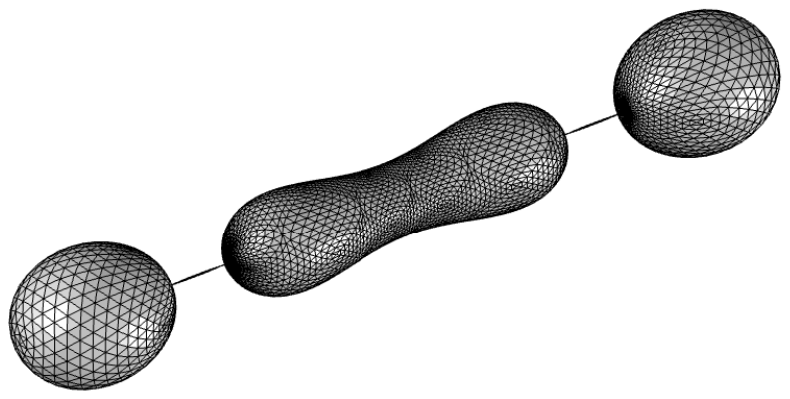}}
	\subfigure[BGN/PC, $t\approx 0.4612$]{\includegraphics[width=.32\textwidth,trim=0 250 0 250, clip]{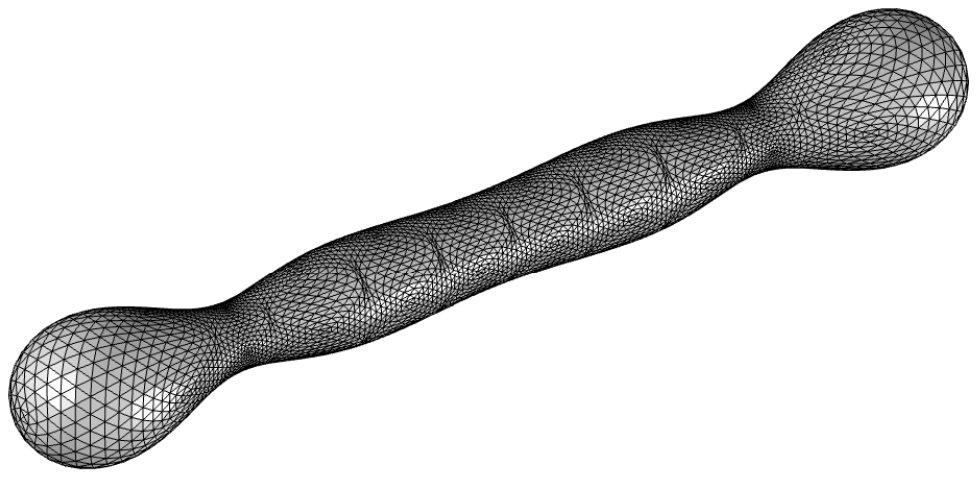}}
	\subfigure[BGN/PC, $t\approx 0.6720$]{\includegraphics[width=.32\textwidth,trim=0 250 0 250, clip]{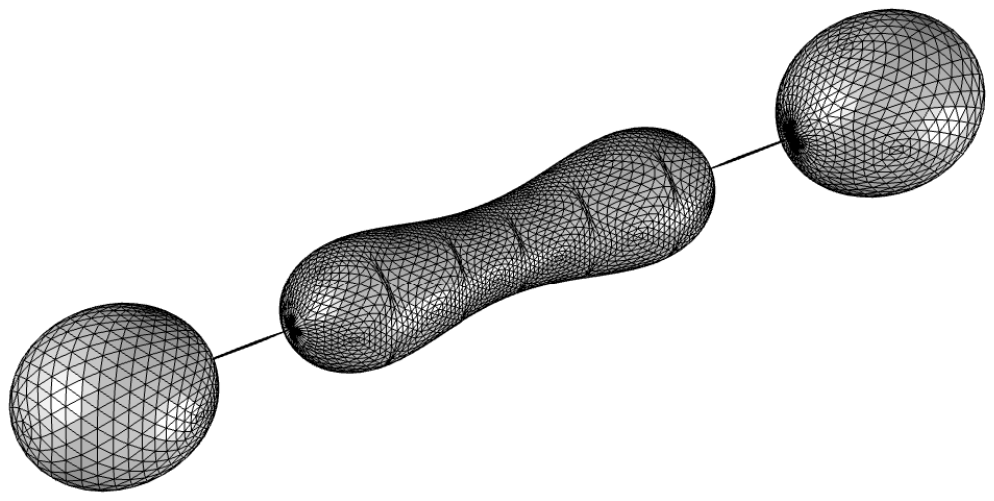}}	
	\caption{Numerical surfaces from different schemes  for \Cref{ex-box}, the images are rescaled.}
	\label{fig-Box}
\end{figure}

\begin{figure}[!htp]
	\centering
	\subfigure[Evolution of surface areas]{\includegraphics[width=.42\textwidth,trim=0 200 0 200, clip]{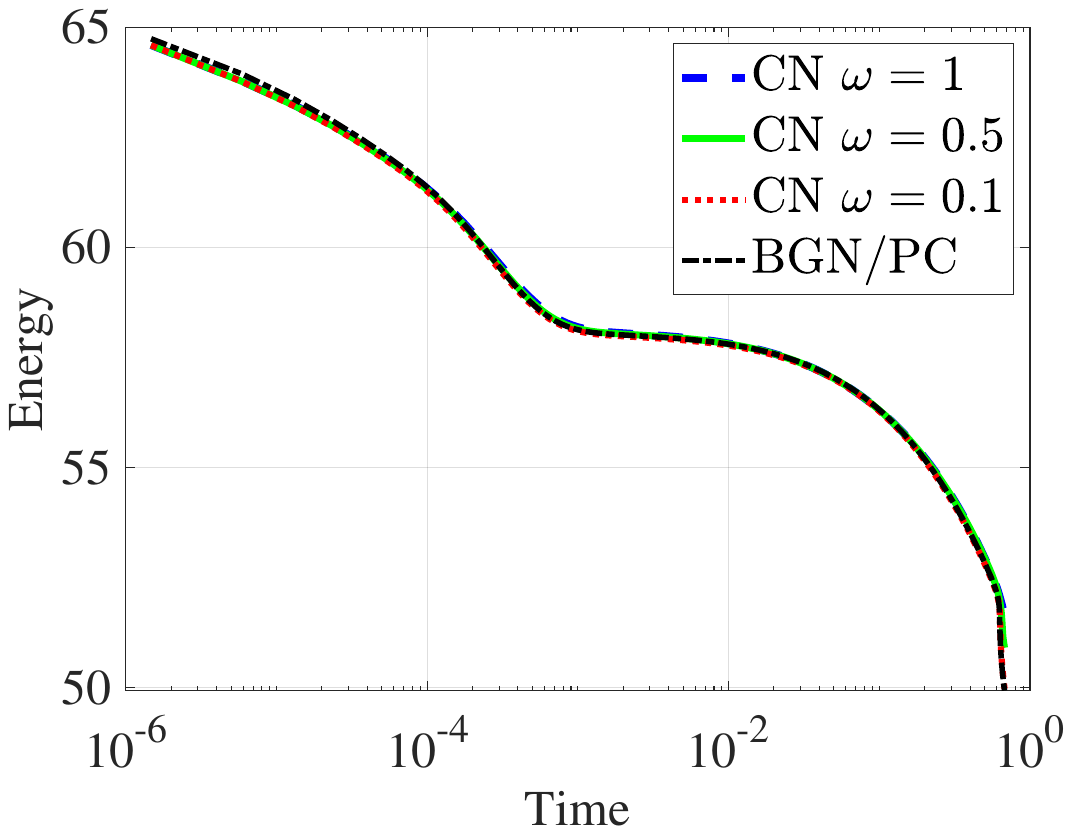}}
	\subfigure[The mesh quality index $\sigma_{\max}$]{\includegraphics[width=.42\textwidth,trim=0 200 0 200, clip]{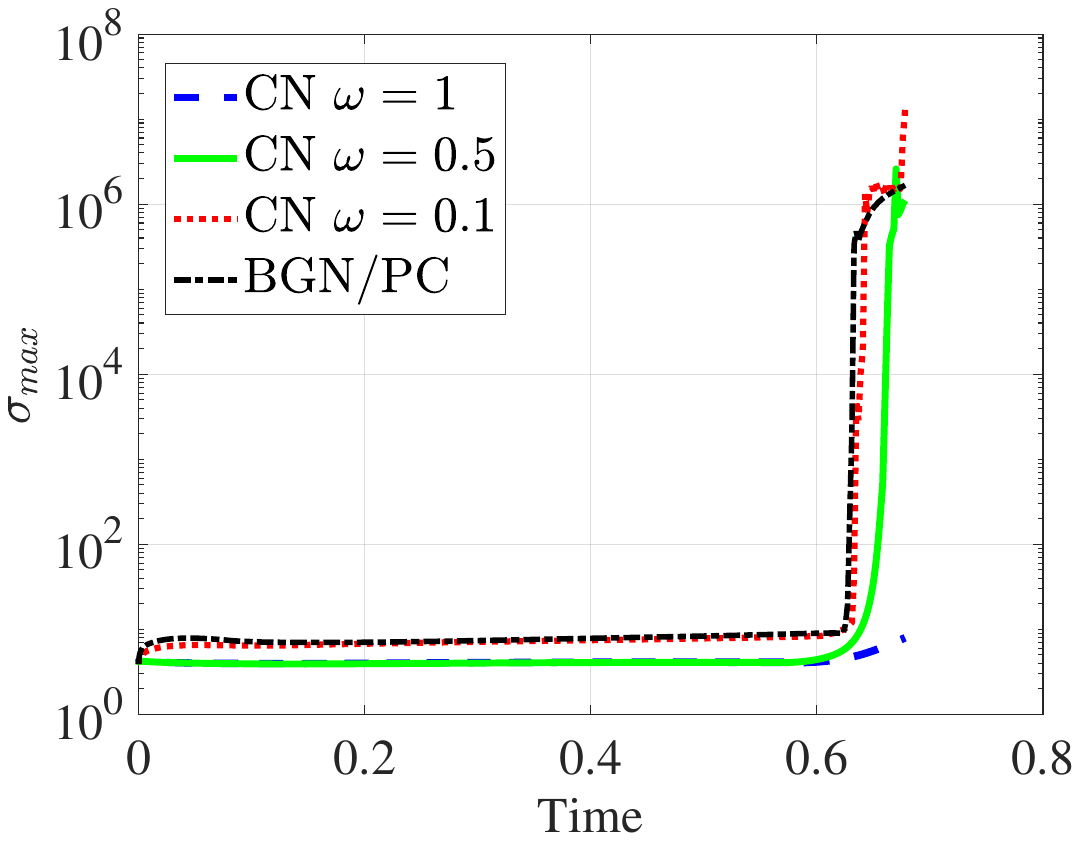}}
	\caption{Surface areas and $\sigma_{\max}$ with respect to time  for \Cref{ex-box}.}
	\label{fig-Box-Energy}
\end{figure}

\begin{example}\label{ex-box}
	In this example, we consider the evolution of a $1:1:16$ box under surface diffusion. The initial mesh, shown in \Cref{fig-Box}(a), consists of $N_p=5378$ vertices and $N_T=10752$ plane triangles. We adopt the graded time mesh $t_m=T(m/M)^2$, $m=0,1,\ldots,M$, with $T=0.68$ and $M=680$, since the initial surface is not smooth. Numerical surfaces at different times obtained from our Crank-Nicolson scheme with $\omega=0.1$ and from the BGN-based predictor--corrector (BGN/PC) scheme proposed in \cite{JIANG2025113901} are presented in \Cref{fig-Box}. One can observe that the mesh of our method remains well distributed before singularities develop, whereas the vertices of the BGN-based scheme cluster undesirably in several places where the surface is smooth. \Cref{fig-Box-Energy}(a) plots the energy curves obtained by the BGN/PC scheme and by our Crank-Nicolson scheme with various values of $\omega$. One can see that they match very well. \Cref{fig-Box}(b) presents $\sigma_{\max}$ for our Crank-Nicolson scheme under different $\omega$ and the BGN/PC scheme. One can see that a bigger $\omega$ can suppress the increase of $\sigma_{\max}$ at the beginning when the initial surface has sharp corners. 
\end{example}


\subsection{Surfaces of positive genus}
In this section we present some results for surfaces of positive genus. As we have pointed out, we can set $\omega=0$ then our Euler scheme reduces to the BGN method and our Crank-Nicolson scheme can be regarded as a modified BGN method which can achieve second-order accuracy in time when the conformal map between $\Gamma(0)$ and $\Gamma(t)$ does exist. 
\begin{example}\label{ex-MCF-g1}
	We first consider the mean curvature with initial surface given by the following parametrization:
	\begin{equation}\label{eq-torus-65}
		\vx=\left[
		\begin{matrix}
			(1+0.65\cos\varphi)\cos\theta \\
			(1+0.65\cos\varphi)\sin\theta \\
			0.65\sin\varphi +0.3\sin 5\theta
		\end{matrix}
		\right],
		\quad \theta\in [0,2\pi) ,\quad\varphi \in [0,2\pi) . 
	\end{equation}
	
	We triangulate the initial surface into $N_T=6752$ plane triangles with $N_p=3376$ vertices, as shown in \Cref{fig-MCF-flower}(a). Numerical surfaces obtained by different schemes with $\tau=1\times 10^{-3}$ at $t=0.235$ are presented in \Cref{fig-MCF-flower}(b)--(e). One can see that the shapes of the four surfaces are very similar. \Cref{fig-MCF-flower-Energy}(a) shows the energy with respect to time. We can see that our methods give almost the same result as the BGN scheme. \Cref{fig-MCF-flower-Energy}(b) plots the mesh quality index $\sigma_{\max}$ for different schemes from which one can see a tiny number of $\omega$ still improve the mesh quality significantly. We mention that  $h_0\approx0.2$ in this example, so  the term involving $\omega$ can be regarded as a perturbation of the equation under $\omega=0$.  
	
	\begin{figure}[!htp]
		\centering
		\subfigure[Initial surface]{\includegraphics[width=.32\textwidth,trim=0 200 0 150, clip]{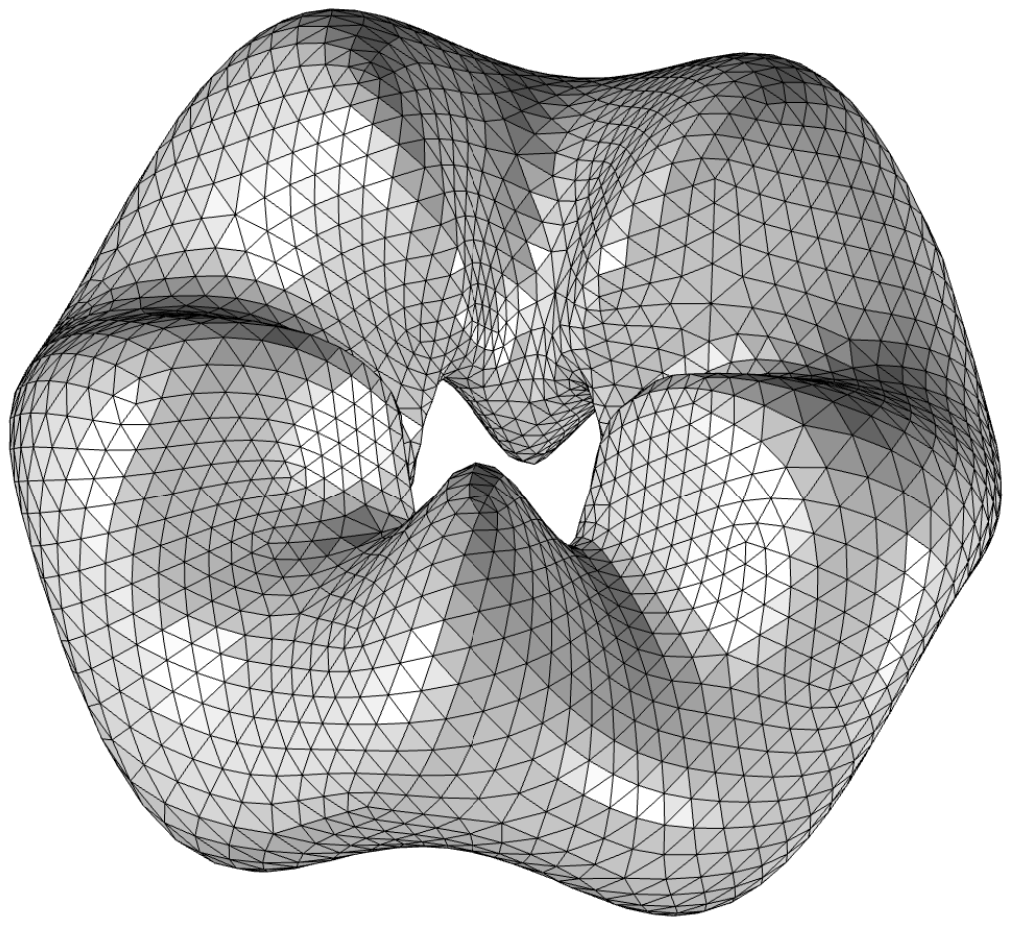}}
		\subfigure[Our CN $\omega=0$]{\includegraphics[width=.32\textwidth,trim=0 200 0 200, clip]{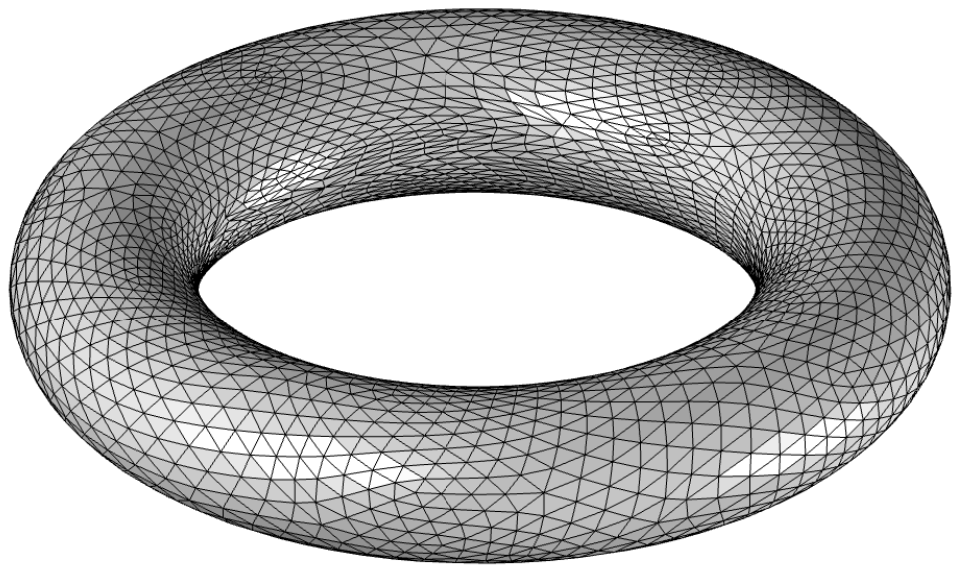}}
		\subfigure[Our CN $\omega=0.01$]{\includegraphics[width=.32\textwidth,trim=0 200 0 200, clip]{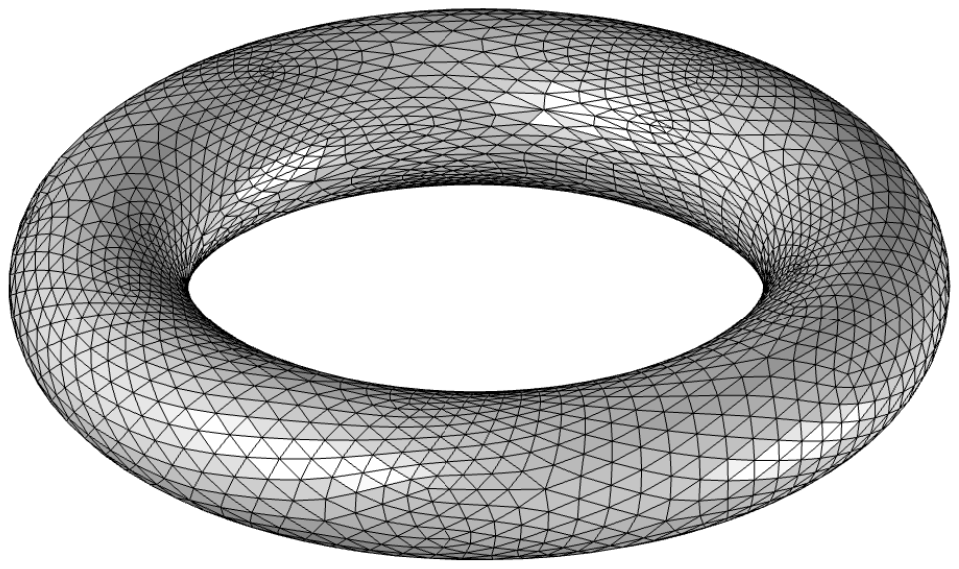}}\\
		\hspace{3.5cm}
		\subfigure[BGN scheme]{\includegraphics[width=.32\textwidth,trim=0 220 0 200, clip]{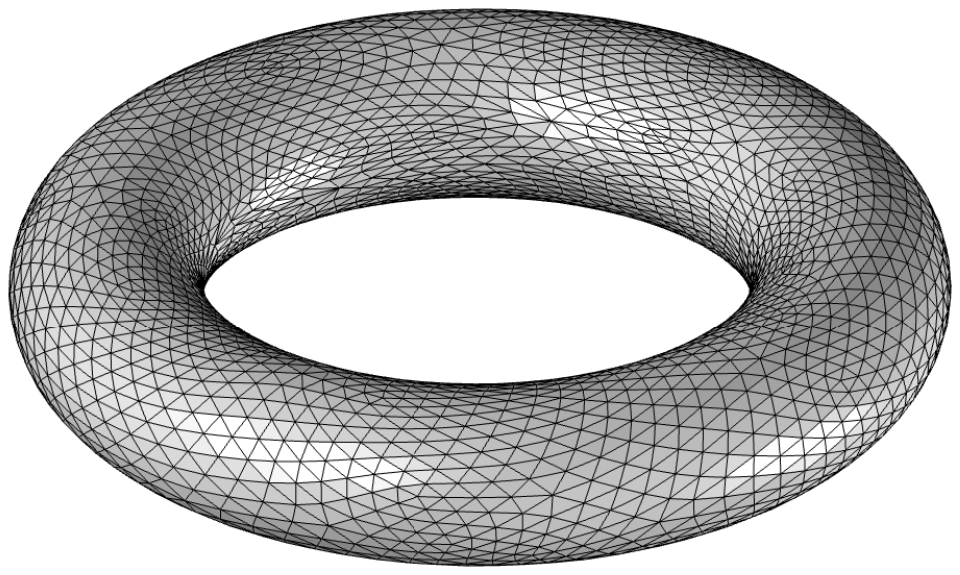}}
		\subfigure[Variant DeTurck method $\alpha^{-1}=2$]{\includegraphics[width=.28\textwidth,trim=0 200 0 200, clip]{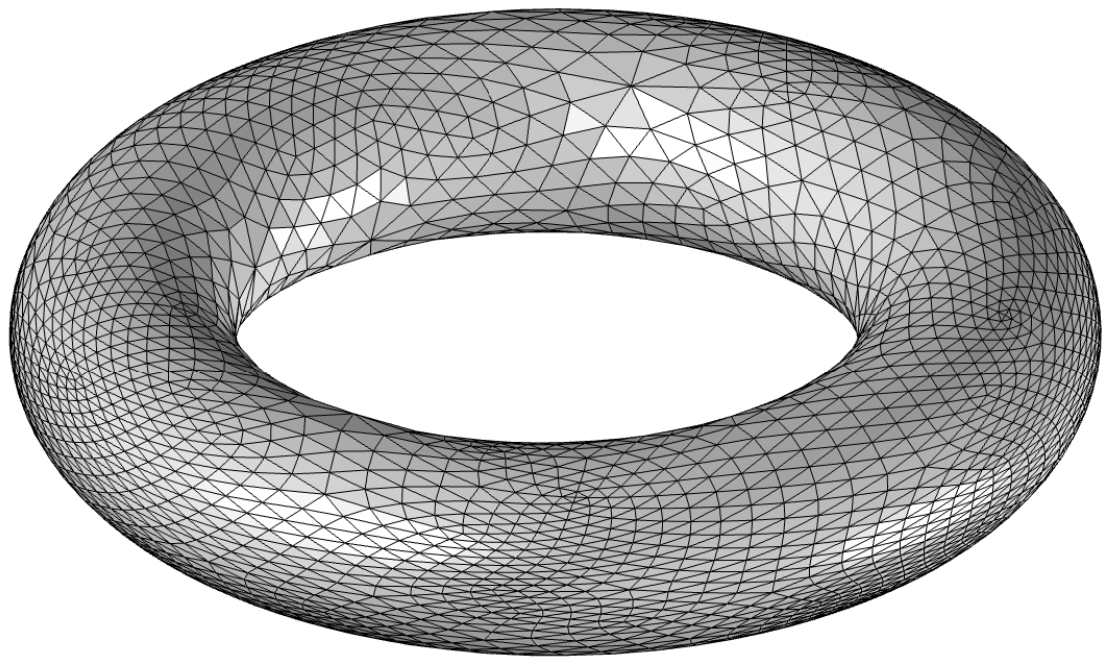}}
		\caption{Numerical surfaces at $t=0.235$ from different schemes with $\tau=1\times 10^{-3} $ for \Cref{ex-MCF-g1}, the images are rescaled. Subplots from (b)--(e): our Crank-Nicolson scheme with $\omega=0$ and $0.01$; the BGN scheme \cite{BGN2008};  the variant DeTurck method with $\cM_h=\Gamma_h^0$, $\alpha^{-1}=2$ \cite[Algorithm 6.4]{Elliott-Fritz-2017}.}
		\label{fig-MCF-flower}
	\end{figure}
	
	\begin{figure}[!htp]
		\centering
		\subfigure[Evolution of surface areas]{\includegraphics[width=.42\textwidth,trim=0 200 0 200, clip]{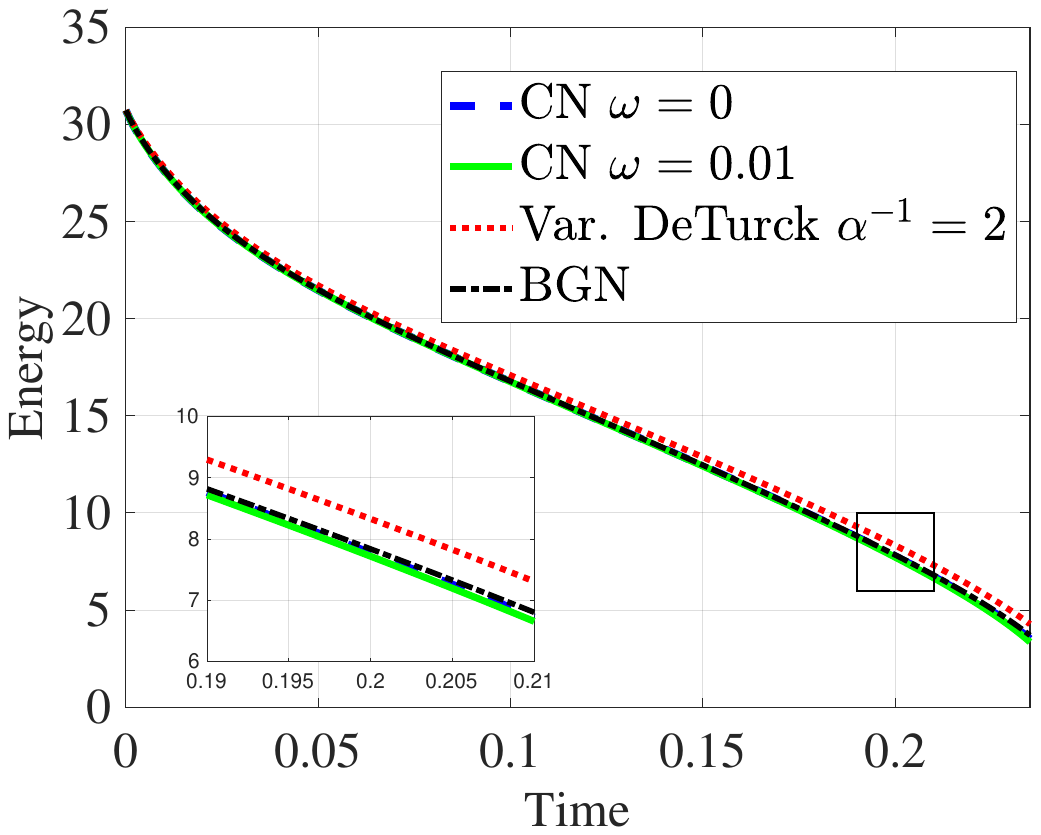}}
		\subfigure[The mesh quality index $\sigma_{\max}$]{\includegraphics[width=.42\textwidth,trim=0 200 0 200, clip]{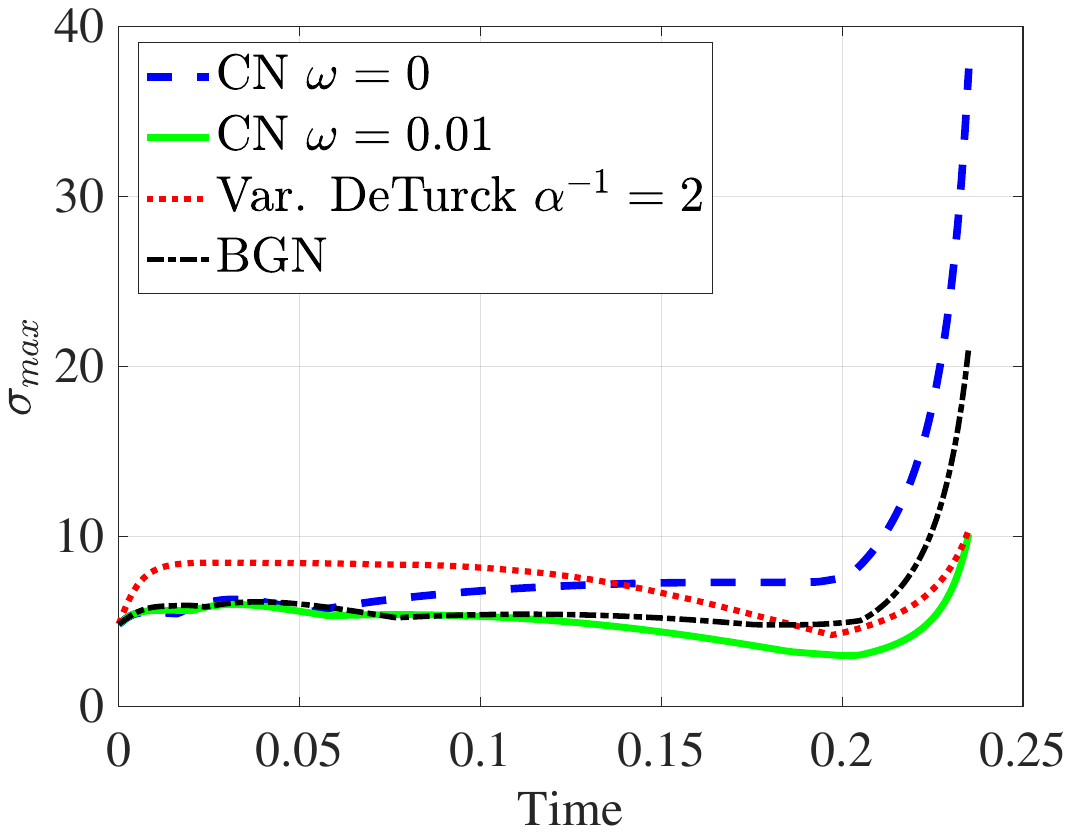}}
		\caption{Surface areas and $\sigma_{\max}$ with respect to time  for \Cref{ex-MCF-g1}. }
		\label{fig-MCF-flower-Energy}
	\end{figure}
\end{example}
\begin{example}\label{ex-genus2}
	In this example, we consider the evolution of a genus-2 surface under surface diffusion. The initial mesh, shown in \Cref{fig-genus-Energy}(a), consists of $N_p=8160$ vertices and $N_T=16324$ plane triangles. We set $\omega=0$ and run our Crank-Nicolson scheme with $\tau=2\times 10^{-4}$ up to $T=0.04$. The numerical surfaces at $t=0.03$ and $t=0.04$ are shown in \Cref{fig-genus2}(b) and (c), respectively. One can observe that the two holes eventually disappear, giving rise to two singularities. \Cref{fig-genus-Energy}(a) presents the energy versus time, and the result of the BGN/PC scheme is also included. One can see that the two results are almost the same before singularities develop, which demonstrates the ability of our method to simulate surfaces beyond topological spheres. \Cref{fig-genus-Energy}(b) compares $\sigma_{\max}$ for our method with that for the BGN/PC method. One can observe that both methods preserve mesh quality well before singularities develop.
	\begin{figure}[!htp]
		\centering
		\subfigure[Initial surface]{\includegraphics[width=.36\textwidth,trim=150 340 150 310, clip]{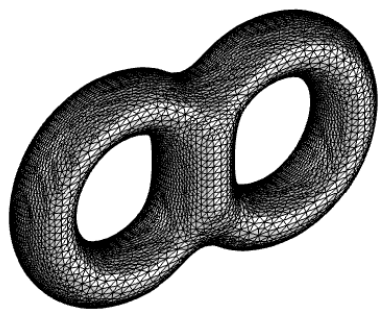}}
		\subfigure[$t= 0.03$]{\includegraphics[width=.31\textwidth,trim=100 280 100 280, clip]{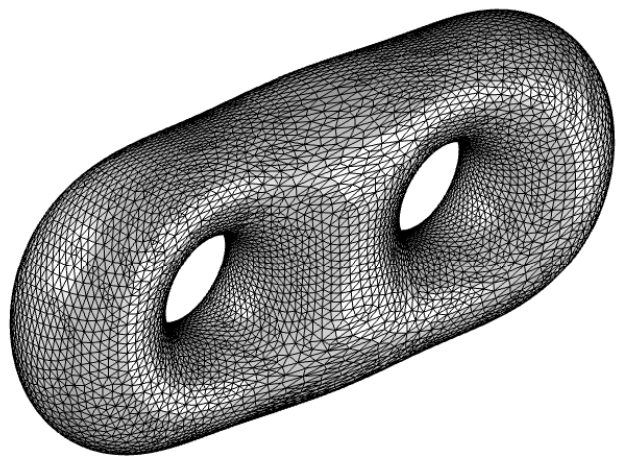}}
		\subfigure[$t= 0.04$]{\includegraphics[width=.31\textwidth,trim=100 280 100 280, clip]{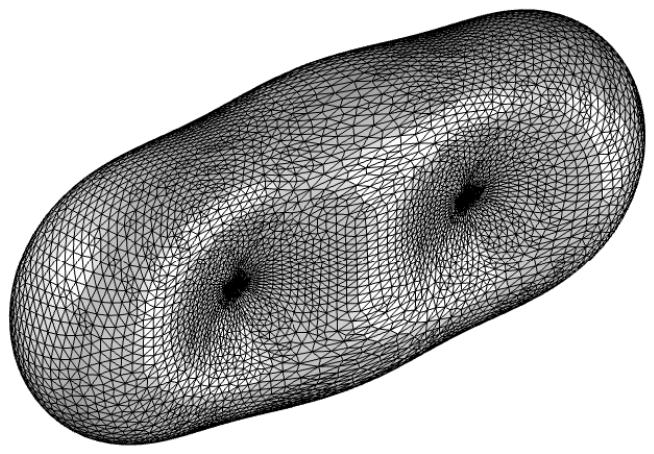}}
		\caption{Numerical surfaces at different times for \Cref{ex-genus2}, the images are rescaled.}
		\label{fig-genus2}
	\end{figure}
	\begin{figure}[!htp]
		\centering
		\subfigure[Evolution of surface areas]{\includegraphics[width=.42\textwidth,trim=0 200 0 200, clip]{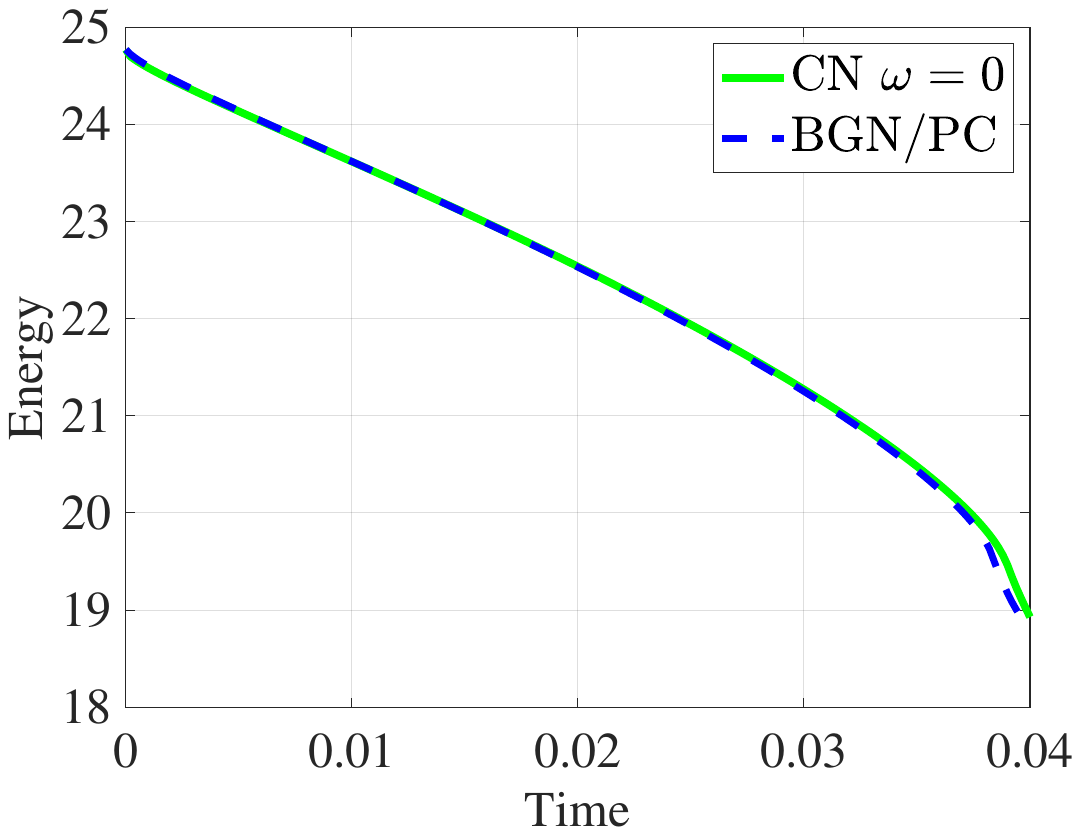}}
		\subfigure[The mesh quality index $\sigma_{\max}$]{\includegraphics[width=.42\textwidth,trim=0 200 0 200, clip]{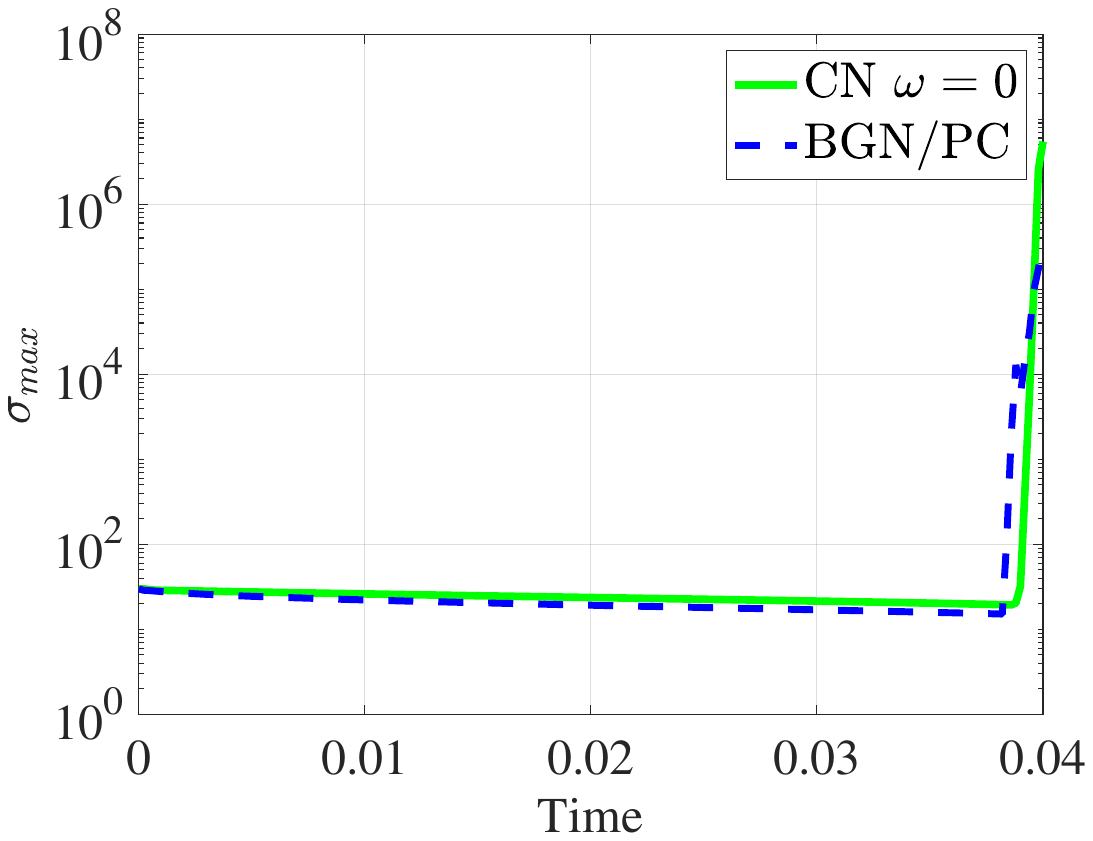}}
		\caption{Surface areas and $\sigma_{\max}$ with respect to time  for \Cref{ex-genus2}.}
		\label{fig-genus-Energy}
	\end{figure}
\end{example}

\section{Conclusion}
We proposed a new method for solving geometric flows of surfaces with or without boundaries. The tangential motion in our method is determined by requiring the flow map to be conformal. We proved the existence and uniqueness of the discrete solutions, as well as the energy stability, for all our schemes. In theory, the conformality of the flow map for surfaces with complicated topology is not guaranteed except for topological spheres and disks.  The advantage of our discrete schemes is that they can be applied to general surfaces, although the second-order scheme may degenerate to first-order accuracy in certain cases. Volume preservation can be easily achieved by using the first- or second-order time-averaged normal vector proposed in \cite{Bao-Zhao-2021} or \cite{duan2026second}.

\bibliographystyle{abbrv}
\bibliography{References_GFs}    
\end{document}